# ASYMPTOTIC EXPANSIONS FOR THE LAPLACE APPROXIMATIONS OF SUMS OF BANACH SPACE-VALUED RANDOM VARIABLES

By Sergio Albeverio and Song Liang[1]

*University of Bonn and Tohoku University*

Let $X_i$, $i \in \mathbf{N}$, be i.i.d. $B$-valued random variables, where $B$ is a real separable Banach space. Let $\Phi$ be a smooth enough mapping from $B$ into $\mathbf{R}$. An asymptotic evaluation of $Z_n = E(\exp(n\Phi(\sum_{i=1}^n X_i/n)))$, up to a factor $(1+o(1))$, has been gotten in Bolthausen [*Probab. Theory Related Fields* **72** (1986) 305–318] and Kusuoka and Liang [*Probab. Theory Related Fields* **116** (2000) 221–238]. In this paper, a detailed asymptotic expansion of $Z_n$ as $n \to \infty$ is given, valid to all orders, and with control on remainders. The results are new even in finite dimensions.

**1. Introduction.** Let $(B, \|\cdot\|_B)$ be a Banach space and $\mu$ be a probability measure on $B$. We assume that the smallest closed affine space that contains $\operatorname{supp}\mu$ is $B$. Moreover, we assume the following:

Assumption A1. There exists a constant $K_1 > 0$ such that

$$\int_B \exp(K_1\|x\|_B^2)\mu(dx) < \infty.$$

(This is satisfied if, e.g., $\mu$ is a Gaussian measure and $K_1 > 0$ is sufficiently small, by Fernique's estimate. See, e.g., [20].)

Let $\Phi : B \to \mathbf{R}$ be a five times continuously Fréchet differentiable function satisfying the following:

Received June 2002; revised September 2004.
[1]Supported in part by Alexander von Humboldt Foundation (Germany) and Grant-in-Aid for the Encouragement of Young Scientists (No. 15740057), Japan Society for the Promotion of Science.
AMS 2000 subject classifications. 62E20, 60F10, 60B12.
*Key words and phrases.* Laplace approximation, asymptotic expansions, i.i.d. random vectors, Banach space-valued random variables.

ASSUMPTION A2. There exists a constant $K_2 > 0$ such that

$$\Phi(x) \le K_2(1 + \|x\|_B) \qquad \text{for all } x \in B.$$

We remark that this is a one-sided condition; it is satisfied, for example, if $\Phi$ is negative for $\|x\|_B \to \infty$.

Let $X_n$ and $S_n$, $n \in \mathbf{N}$, be the random variables defined by $X_n(\underline{x}) = x_n$ and $S_n(\underline{x}) = \sum_{k=1}^n x_k$ for any $\underline{x} = (x_1, x_2, x_3, \dots) \in B^{\mathbf{N}}$.

We are interested in the behavior of

$$Z_n \equiv E^{\mu^{\otimes\infty}}\left[\exp\left(n\Phi\left(\frac{S_n}{n}\right)\right)\right] \qquad \text{as } n \to \infty,$$

where $E^\lambda$ stands for the expectation with respect to the measure $\lambda$, and $\mu^{\otimes\infty}$ is the product of $|\mathbf{N}|$ copies of $\mu$ [corresponding to the distribution of the $(X_n)_{n\in\mathbf{N}}$].

By Donsker and Varadhan [13], we have that

$$\lim_{n\to\infty} \frac{1}{n} \log Z_n = \sup_{x\in B}\{\Phi(x) - h(x)\},$$

where $h$ is the entropy function of $\mu$:

$$h(x) = \sup_{\phi\in B^*}\{\phi(x) - \log M(\phi)\}, \qquad x \in B,$$

$B^*$ is the dual space of $B$ and $M(\phi) = \int_B e^{\phi(x)}\mu(dx)$ for any $\phi \in B^*$.

It has been shown by Bolthausen [8] that there is at least one $x^* \in B$ with $\Phi(x^*) - h(x^*) = \sup_{x\in B}\{\Phi(x) - h(x)\}$, and the set $K = \{x \in B; \Phi(x) - h(x) = \sup_B\{\Phi - h\}\}$ is compact. Also, we assume the following, as in [8] and [21].

ASSUMPTION A3. There is a unique $x^* \in B$ with $\Phi(x^*) - h(x^*) = \sup_{x\in B}\{\Phi(x) - h(x)\}$.

This is satisfied, for example, if $\Phi$ is strictly concave, since $h$ is always convex by the definition of it.

We will use $x^*$ exclusively for this point.

Let $\nu$ be the probability measure on $B$ given by

$$\nu(dx) = \frac{\exp(D\Phi(x^*)(x))\mu(dx)}{M(D\Phi(x^*))},$$

where $D$ means the Fréchet derivative. As has been shown by Bolthausen [8], the following proposition holds.

PROPOSITION 1.1. *Under Assumptions* A1–A3,

$$x^* = \int_B x\nu(dx), \tag{1.1}$$

$$h(x^*) = D\Phi(x^*)(x^*) - \log M(D\Phi(x^*)). \tag{1.2}$$

Let $\nu_0$ be the 0-centered measure associated with $\nu$, that is, $\nu_0 = \nu\theta_{x^*}^{-1}$, where $\theta_a : B \to B$ is defined by $\theta_a(x) = x - a$, $x \in B$.

Let $\Gamma(\varphi, \psi) = \int_B \varphi(x)\psi(x)\nu_0(dx)$ be the covariance (under $\nu_0$) of $\varphi$ and $\psi$ for any $\varphi, \psi \in B^*$. Then $\Gamma$ becomes an inner product on $B^*$. Let $H \equiv (\overline{B^*}^{\Gamma})^*$, where $\overline{B^*}^{\Gamma}$ means the completion of $B^*$ with respect to $\Gamma$. (It has been shown in [21] that $H$ can be regarded as a dense subset of $B$.)

The following holds, as shown in [8]:

$$D^2\Phi(x^*)(\iota(\phi), \iota(\phi)) \le \Gamma(\phi, \phi) \qquad \text{for any } \phi \in B^*,$$

where $\iota(\phi) \equiv \int_B \phi(x)x\nu_0(dx)$, $\phi \in B^*$. From this, we see that all of the eigenvalues of the symmetric operator $D^2\Phi(x^*)|_{H\times H}$ [given by the restriction of $D^2\Phi(x^*)$ to $H \times H$] are not greater than 1. We assume the following as in [8] and [21].

ASSUMPTION A4. All of the eigenvalues of $D^2\Phi(x^*)|_{H\times H}$ are strictly smaller than 1.

This is a nondegeneracy condition (depending on both $\Phi$ and $\mu$), which says that $\Phi - h$ has a nonvanishing curvature at the point $x^*$ (see [8]), or, by the words of [23], means that the Hessian of $\Phi$ is strictly positive definite. This condition also implies that the determinant appearing in the following Proposition 1.2 is different from 0. In the sense of the theory of singularities of maps (see, e.g., [5, 26]) this is a genericity condition (i.e., if it were not satisfied for a given $\Phi$, a "small perturbation" of $\Phi$ would make it satisfied; see, e.g., [24]).

Bolthausen [8] and Kusuoka and Liang [21] studied the leading term of $\exp(-n(\Phi(x^*) - h(x^*)))Z_n$ as $n \to \infty$. In particular, [21] gives us the following.

PROPOSITION 1.2. *Let Assumptions* A1–A4 *be satisfied, and assume some technical condition for controlling the third remainder of* $\Phi$ *in the Taylor expansion around* $x^*$ *(see* [21], (A5) *for the explicit expression of this condition, which will be replaced by a stronger Assumption* A5*). Then we have*

$$\lim_{n\to\infty} \exp(-n(\Phi(x^*) - h(x^*)))Z_n$$
$$= \exp\left(\int_B D^2\Phi(x^*)(y,y)\nu_0(dy)\right) \det{}_2(I_H - D^2\Phi(x^*))^{-1/2} \equiv C_0(x^*).$$

Note that both [21] and our present paper do not assume the so-called "Central Limit Theorem Assumption" used in [8], which restricts the spaces.

Now it is a natural problem to investigate the more precise asymptotic behavior of $Z_n$ as $n \to \infty$, beyond the leading term. This does not seem to have been discussed before and is entirely in the spirit of corresponding investigations for the case of real-valued random variables, where one wishes to go beyond the functional and central limit theorem, that is, in the sense of Edgeworth expansions, for a certain functional of the normalized sum variables $S_n/\sqrt{n}$; see, for example, [7, 18, 19] and references therein. Even in the case of real random variables, our results, however, are not reduced to known results, because of the form of our function $\phi$; see Remark 1.4 below. The aim of this paper is to answer this question in our general setting of Assumptions A1–A4, adding Assumptions A5 and A6 (which just are a little stronger than what follows from Assumptions A1–A4).

We state now the condition implying the one from [21] mentioned in Proposition 1.2.

Assumption A5. There exists a constant $\delta_0 > 0$ and a bilinear, symmetric, bounded function $K_5 : B \times B \to \mathbf{R}$ such that

$$|D^5\Phi(x)(y,y,y,y,y)| \le \|y\|_B^3 K_5(y,y)$$

for any $x \in B$ with $\|x - x^*\|_B < \delta_0$ and any $y \in B$.

We remark that this implies the technical condition A5 in [21] about the third remainder of the Taylor expansion of $\Phi$, so that in particular, we can use Proposition 1.2 under the sole Assumptions (A1)–(A5).

Assumption A5 is satisfied if, for example, $B$ is a separable real Hilbert space. Actually, if $B$ is a real Hilbert space, writing the inner product of $B$ as $(\cdot,\cdot)_B$, then we can just take $K_5(x,y) = C(x,y)_B$ for $x, y \in B$, with $C$ the supremum of the operator norm of $D^5\Phi(x)$ on $\{x \in B; \|x - x^*\| \le \delta_0\}$.

For the sake of simplicity, we denote $D^i\Phi(x^*)$ by $\Psi_i$, $i = 2,3,4$. We have by [21] that $\Psi_2|_{H\times H}$ is a Hilbert–Schmidt operator, hence the corresponding resolvent set consists only of eigenvalues. We shall denote the eigenvalues by $a_k$, $k \in \mathbf{N}$, and the corresponding eigenfunctions by $e_k, k \in \mathbf{N}$. Without loss of generality, we may assume that $e_k, k \in \mathbf{N}$, consists of an orthonormal base (ONB) of the dual $H^*$ of $H$. Let $f_k, k \in \mathbf{N}$, be the corresponding ONB of $H$. Also, by [21] (for those $k \in \mathbf{N}$ with $a_k \ne 0$), we may assume that $e_k \in B^*$.

Then we have

$$\Psi_2(x,y) = \sum_{k=1}^{\infty} a_k e_k(x) e_k(y) \tag{1.3}$$

for any $x, y \in H$. We remark at this place that by Minlos' theorem, the same equation holds for $\nu_0$-a.s. $x, y \in B$ if $\Psi_2|_{H\times H}$ is a nuclear operator, that is, if $\sum_{k=1}^{\infty} |a_k| < \infty$.

Now, we are able to formulate our last assumption.

ASSUMPTION A6. There exists a bilinear, symmetric, bounded function $\widetilde{\Psi}_2 : B \times B \to \mathbf{R}$, and a monotone nonincreasing sequence of positive numbers $\delta_N, N \in \mathbf{N}$, that converges to 0 as $N \to \infty$ such that for any $N \in \mathbf{N}$,

$$\sum_{k>N,\, a_k>0} a_k e_k^{\otimes 2} \le \delta_N \widetilde{\Psi}_2.$$

REMARK 1.1. Assumption A6 implies that $\sum_{k=1}^{\infty} |a_k| e_k \otimes e_k$ is well defined and gives a continuous operator on $B \times B$. (Actually, in this case, we have that both $\sum_{k=1}^{\infty} a_k e_k^{\otimes 2}|_{H\times H}$ and $\sum_{a_k>0} a_k e_k^{\otimes 2}|_{H\times H}$, and hence also $\sum_{k=1}^{\infty} |a_k| e_k^{\otimes 2}|_{H\times H}$, are continuous with respect to the $B$-norm, so we can extend them in a continuous way to the whole $B$.) Therefore, $\sum_{k=1}^{\infty} |a_k| < \infty$, that is, $\Psi_2|_{H\times H}$ is a nuclear operator, and hence as already remarked above, (1.3) holds for $\nu_0$-a.s. $x, y \in B$.

We emphasize here that $\{a_k\}_{k\in\mathbf{N}}$ and $\{e_k\}_{k\in\mathbf{N}}$ are eigenvalues and the corresponding eigenfunctions of $\Psi_2$ acting in $H$, instead of $B$. And $B$ and $H$ are different even if $B$ is a Hilbert space.

We remark that Assumption A6 is satisfied, for example, if in the representation (1.3) of $\Psi_2 = D^2\Phi(x^*)$, there is only a finite number of $a_k$ which are strictly positive, or more generally, if there exists a $p > 1$ such that $\sum_{k\,:\,a_k>0} a_k^{1/p} e_k^{\otimes 2}$ is a bounded function on $B \times B$. In other words, for any bounded positive-definite function $A : B \times B \to \mathbf{R}$, write the eigenvalues and corresponding eigenfunctions as $b_k$ and $\tilde{e}_k$, $k \in \mathbf{N}$, that is, $A = \sum_{k=1}^{\infty} b_k \tilde{e}_k^{\otimes 2}$; then it is easy to see that $b_k \to 0$ as $k \to \infty$, so for any $p > 1$, we have that $\sum_{k=1}^{\infty} b_k^p \tilde{e}_k^{\otimes 2}$ satisfies our Assumption A6. See also the discussion following Theorem 1.3.

Also, see the end of this section for examples where all Assumptions A1–A6 are satisfied.

As in [21], let $H_1$ be a Hilbert space that includes $H$ as a subset with the embedding being a Hilbert–Schmidt operator. (See Section 3 for the precise definition and the construction of $H_1$ using Assumption A6.) Then there exists an $H_1$-valued Gaussian random variable $Y$ such that the distribution of $u(Y)$ is $N(0, \|u\|_{H^*}^2)$ (the normal distribution with mean 0 and variance $\|u\|_{H^*}^2$) for any $u \in H^*$. Since

$$E^Y\left[E^{\nu_0}\left[\left|\sum_{k=1}^{\infty} \sqrt{a_k} e_k(X_1) e_k(Y)\right|^2\right]\right] = E^{\nu_0}\left[\sum_{k=1}^{\infty} |a_k| e_k(X_1)^2\right]$$

$$= \sum_{k=1}^{\infty} |a_k| < \infty$$

(with $E^Y$ the expectation with respect to the distribution $P_Y$ of $Y$), we have that $|\sum_{k=1}^{\infty} \sqrt{a_k} e_k(X_1) e_k(Y)| < \infty$ for a.e.-$(X_1, Y)$ with respect to the measure $\mu \otimes P_Y$. We shall write $\sum_{k=1}^{\infty} \sqrt{a_k} e_k(X_1) e_k(Y)$ as $(\Psi_2^{1/2} X_1, Y)$. We remark that this may be a complex number since the coefficients $a_k$ may be negative. Also, it is easy to see by the central limit theorem in Hilbert spaces (e.g., [4]) that $\frac{S_n}{\sqrt{n}} \to Y$ in law under $\nu_0^{\otimes \infty}$ as $n \to \infty$, and since $\Psi_2|_{H \times H}$ is nuclear under Assumption A6, the constant $C_0(x^*)$ in Proposition 1.2 is equal to $\det(I_H - \Psi_2)^{-1/2}$ (where det is the Fredholm determinant).

Now, we are ready to give our main result, which provides a precise expression for the coefficient of the term $n^{-1}$ in the expansion of $U_n := \exp(-n(\Phi(x^*) - h(x^*)))Z_n$.

THEOREM 1.3. *Under Assumptions A1–A6 above, we have that*

$$\lim_{n \to \infty} n(\exp(-n(\Phi(x^*) - h(x^*)))Z_n - \det(I_H - \Psi_2)^{-1/2}) = C_2(x^*),$$

*where* $C_2(x^*)$ *is the constant given by*

$$\begin{aligned}
C_2(x^*) = {} & E^Y\Big[e^{\Psi_2(Y,Y)/2}\Big(-\frac{1}{8}\Psi_2(Y,Y)^2 + \frac{1}{2(3!)^2} E^{\nu_0}[(\Psi_2^{1/2} X_1, Y)^3]^2 \\
& \qquad + \frac{1}{4!} E^{\nu_0}[(\Psi_2^{1/2} X_1, Y)^4]\Big)\Big] \\
& + E\Big[e^{\Psi_2(Y,Y)/2}\Big(\frac{1}{2(3!)^2}\Psi_3(Y,Y,Y)^2 + \frac{1}{4!}\Psi_4(Y,Y,Y,Y)\Big)\Big] \\
& + \frac{1}{3!} E^Y[e^{\Psi_2(Y,Y)/2}] E^{\nu_0}[\Psi_3(X_1, X_1, X_1)] \\
& + \sum_{k=3}^{4} \frac{1}{(k-1)!} \\
& \quad \times \sum_{i_1 + \cdots + i_{k-1} + i_k/3 = 2(k-2)} \Big(\frac{1}{2}\Big)^{\sum_{j=1}^{k-1}(i_j - 1)} \Big(\frac{1}{3!}\Big)^{i_k/3} \\
& \qquad \times E^Y\Big[e^{\Psi_2(Y,Y)/2} E^{\nu_0^{\otimes k}} \\
& \qquad\quad \times \Big[\prod_{j=1}^{k} (\Psi_2^{1/2} X_j, Y)^{i_j} \Psi_3(X_1, X_2, X_{k-1})\Big]\Big].
\end{aligned}$$

REMARK 1.2. Ellis and Rosen [17] considered the same problem of "large $n$ expansion," but only for the Gaussian case, that is, when $\mu$ is a Gaussian measure on some functional space (e.g., $L^2$-space), and their method used the fact of having Gaussian measures in an essential way. (This work continues previous works on Laplace method for infinite-dimensional Gaussian measures by Pincus, Schilder and Donsker and Varadhan. See references in [17] and, e.g., [3, 2].)

In the special case described in [17], that is, the Gaussian case, our Assumption A6 can be rewritten as follows: Let $A$ denote the covariance of the Gaussian measure $\mu$ on Banach space $B$; then our $\nu_0$ is nothing but $N(0,A)$. Consider $D^2\Phi(x^*)(A^{1/2}\cdot, A^{1/2}\cdot) : B \times B \to \mathbf{R}$. It is easy to see that this is a nuclear operator. Let $\{\tilde{a}_k\}_{k\in\mathbf{N}}$ and $\{u_k\}_{k\in\mathbf{N}} \subset B^*$ be the eigenvalues and the corresponding eigenfunctions of it. Then as before, without loss of generality, we may assume that $\{u_k\}_{k\in\mathbf{N}}$ consists of an ONB of $B^*$. Also, $A^{-1/2}u_k \in H^*$, and by extension if necessary, we may assume that $A^{-1/2}u_k \in B^*$. Now, since $\{\tilde{a}_k\}_{k\in\mathbf{N}}$ and $\{A^{-1/2}u_k\}_{k\in\mathbf{N}}$ are the eigenvalues and the corresponding eigenfunctions of $D^2\Phi(x^*)|_{H\times H}$, we have that our Assumption A6 can be implied by the following condition: there exists a $p > 1$ such that $\sum_{k\in\mathbf{N}:\,\tilde{a}_k>0} \tilde{a}_k^{1/p}(A^{-1/2}u_k)^{\otimes 2}$ is a continuous function on $B \times B$.

We stress that our Theorem 1.3 holds without any assumption on $\mu$ to be a Gaussian measure.

The basic idea of our proof is to use the fact that the Laplace transform of a Gaussian measure is exp(quadratic form). With the help of this observation, we then use the independence of $X_k, k \in \mathbf{N}$, to discuss the a.s.-convergence and the dominations, which then implies the $L^1$-convergence, and hence the stated asymptotic formula.

REMARK 1.3. As remarked before the statement of Theorem 1.3, this theorem gives the coefficient of the term $n^{-1}$ in the expansion of $U_n = \exp(-n(\Phi(x^*) - h(x^*)))Z_n$ in powers of $(\frac{1}{n})^{1/2}$. The same kind of result is not known, to the best of the authors' knowledge, even for the finite-dimensional case. By using the same method, we can also give the explicit expression of the coefficient $C_N(x^*)$ of the term $n^{-N/2}$ for any $N \geq 2$ in the expansion of $U_n$, under natural assumptions about the smoothness of $\Phi$ and an assumption corresponding to Assumption A5. We do not write this explicit expansion in this paper, because of its complicated form, but have limited ourselves to explaining our method, taking the case of $C_2(x^*)$ as an example. We rather limit ourselves to give, in Section 3 (cf. Theorem 3.14), the expansion, to all orders in $n^{-N/2}$, of the term $E^{\nu_0^{\otimes\infty}}[\exp(\frac{n}{2}\Psi_2(\frac{S_n}{n}, \frac{S_n}{n})), \|\frac{S_n}{n}\|_B < \varepsilon]$, for $\varepsilon > 0$ small enough. See Remark 1.4 and Section 3.

REMARK 1.4. In this paper we also give, in particular, the asymptotic expansion of $E^{\nu_0^{\otimes\infty}}[\exp(\frac{n}{2}\Psi_2(\frac{S_n}{n},\frac{S_n}{n})), \|\frac{S_n}{n}\|_B<\varepsilon]$, for $\varepsilon>0$ small enough, to any order, with controls on remainders (cf. Theorem 3.14 and Remark 3.1).

In the sense explained in Remark 1.3, we have got an analogue of the Edgeworth expansion for the functional $\exp(\frac{1}{2}\Psi_2)$, with $\Psi_2$ a bilinear, symmetric and bounded function on $B\times B$ that satisfies Assumptions A4 and A6, of the normalized sum variables $S_n/\sqrt{n}$. The Edgeworth expansion with respect to the distribution function in the finite-dimensional case $B=\mathbf{R}$ has been obtained by many authors (see, e.g., [7, 18, 19] and references therein), but all of these give only estimations which are uniform with respect to the variable of the distribution function, and are not usable in the case of our problem (because of the lack of integrability of the function exp with respect to the Lebesgue measure). In fact, it is easy to see that, for example, in the expression

$$\begin{aligned} n\biggl(E^{\nu_0^{\otimes\infty}}\biggl[\exp\biggl(\frac{1}{2}\Psi_2\biggl(\frac{S_n}{\sqrt{n}},\frac{S_n}{\sqrt{n}}\biggr)\biggr),\biggl|\Psi_2\biggl(\frac{S_n}{\sqrt{n}},\frac{S_n}{\sqrt{n}}\biggr)\biggr|\le n\varepsilon\biggr] \\ - E\biggl[\exp\biggl(\frac{1}{2}\Psi_2(Y,Y)\biggr)\biggr]\biggr) \\ = \frac{1}{2}\int_{y\in[-n\varepsilon,n\varepsilon]} e^{y/2} n\biggl(P^{\nu_0^{\otimes\infty}}\biggl(y\le\Psi_2\biggl(\frac{S_n}{\sqrt{n}},\frac{S_n}{\sqrt{n}}\biggr)\le n\varepsilon\biggr) \\ - P(y\le\Psi_2(Y,Y)\le n\varepsilon)\biggr)\,dy \\ + ne^{-n\varepsilon/2}P^{\nu_0^{\otimes\infty}}\biggl(\biggl|\Psi_2\biggl(\frac{S_n}{\sqrt{n}},\frac{S_n}{\sqrt{n}}\biggr)\biggr|\le n\varepsilon\biggr) - ne^{-n\varepsilon/2}P(|\Psi_2(Y,Y)|\le n\varepsilon) \\ - nE\biggl[\exp\biggl(\frac{1}{2}\Psi_2(Y,Y)\biggr),|\Psi_2(Y,Y)|\ge n\varepsilon\biggr], \end{aligned}$$

all the terms except the first one on the right-hand side decay exponentially as $n\to\infty$; hence a uniform estimation of $n(P^{\nu_0^{\otimes\infty}}(y\le\Psi_2(\frac{S_n}{\sqrt{n}},\frac{S_n}{\sqrt{n}})\le n\varepsilon) - P(y\le\Psi_2(Y,Y)\le n\varepsilon))$ with respect to $y\in[-n\varepsilon,n\varepsilon]$ is not enough to obtain the asymptotic expansion we give in Section 3.

REMARK 1.5. In this paper we concentrate on providing asymptotic expansions for the nondegenerate case (much in the spirit of corresponding investigations using Laplace method for functionals of Brownian motion, see, e.g., [6]; it should, however, be stressed that we concentrate on limits of sums of random variables, not on the limit of their distributions). We plan to extend the results to the degenerate case in subsequent publications (for first results on the limit theorem in this case, see [9] and [24]). The investigation of this paper belongs to the general area of probability theory which investigates the asymptotics of processes. See, for example, [11, 14, 22, 25] for the connection with questions of asymptotics for continuous-time processes. The latter two papers deal in particular with the leading term of a Laplace approximation of diffusions, and **(year?)** includes Brownian motion on tori. Expansions beyond the leading term in these "continuum cases," in the generality of our present paper, have not yet been obtained; our paper can also be seen as a first step in this direction.

There are also relations in motivations and some of the methods with other works on asymptotics; see, for example, [[1, 2, 10, 27]], and references therein.

Finally, let us illustrate the use of our main result in a model of classical statistical mechanics. Consider a system of $n$ particles, with the distribution of the state of each particle being given by a probability measure $\mu_0$ on a compact set $M$. Suppose that the interaction of the two particles with states $x$, respectively $y$, is $\frac{1}{n}V(x,y)$, $x,y\in M$, for some "nice" real-valued function $V$ on $M\times M$, and for given $n\in\mathbf{N}$. Then the probability of the system to be in the state given by a Borel subset $A$ of $M^n$ is $\int_A \nu_n(d\underline{x}) = Z_n^{-1}\int_A e^{(1/n)\sum_{i,j=1}^n V(x_i,x_j)}\mu_0^{\otimes n}(dx_1\cdots dx_n)$, $Z_n$ being the normalizing constant [$\underline{x}$ stands for $(x_1,\dots,x_n)\in M^n$]. Relevant interesting physical quantities can be expressed as expectations of the form $E^{\nu_n}[\sum_{i_1,\dots,i_m=1}^n f(X_{i_1},\dots,X_{i_m})]$, for some bounded continuous "observable" function $f:M^m\to\mathbf{R}$. The problem of computation of such expectations as $n\to\infty$ can be generalized to the following one. Let $B$ be equal to the topological dual $C(M)^*$ of $C(M)$, let $X_i=\delta_{x_i}$ and let $\mu$ be the image of $\mu_0$ under $\delta_x$ (looked upon as an element in $B$). Set $\Phi(R)\equiv\int\int V(x,y)R(dx)R(dy)$, $F(R)\equiv\int\cdots\int f(x_1,\dots,x_m)R(dx_1)\cdots R(dx_m)$, $R$ being a positive measure on $M$. Then the above problem can be seen as a particular case of the study of expectations of the form

$$\frac{E^{\mu^{\otimes\infty}}[F((1/n)\sum_{i=1}^n X_i)e^{n\Phi((1/n)\sum_{i=1}^n X_i)}]}{E^{\mu^{\otimes\infty}}[e^{n\Phi((1/n)\sum_{i=1}^n X_i)}]}$$

as $n\to\infty$, where $\mu$ is a probability measure on some Banach space, $F$, $\Phi$ are "good" functions on $B$ and $X_i$ are i.i.d. random variables with distribution $\mu$ on $B$. Since the method for the numerator is exactly the same as that for the denominator (for $F$ smooth), just with the expression more complicated, we limit ourselves to the study of the denominator.

Let us give some more concrete example that satisfies all of our conditions. In the example just given, let $M=\mathbf{T}(=\mathbf{R}/2\pi\mathbf{Z})$, let $\mu_0$ be the uniform distribution on $\mathbf{T}$, $\mu_0(dy)=\frac{1}{2\pi}dy$, and let $V(x,y)=CU(x-y)$, with $U$ a continuous function on $\mathbf{T}$, and $C$ a constant such that $\int_0^{2\pi}\int_0^{2\pi}V(x,y)^2\,dx\,dy\le\pi^2$. Then it is trivial that $D^3\Phi=0$, so our Assumption A5 is satisfied trivially.

Also, the corresponding entropy function is the relative entropy with respect to $\mu_0$ given by

$$h(R) = \int \left(\log \frac{dR}{d\mu_0}\right) dR.$$

So by calculation, we have that the uniform distribution on $\mathbf{T}$ maximizes $\Phi - h$, so the eigenvalues of $D^2\Phi(x^*)$ are nothing but some (global) constant times the coefficients of the Fourier expansion of $U$. Therefore, Assumption A4 is satisfied if $C$ is small enough, and Assumption A6 is also satisfied whenever the Fourier coefficients of $U$ are in $\ell^\alpha$ for some $\alpha \in (0,1)$. Hence in this case all Assumptions A1–A6 are satisfied, and so our theorem applies. Let us remark that this example is related to the mean field model studied in, for example, [15], but for the physically particularly interesting case of translation-invariant interactions. It can also be considered as an generalization of the continuous spinning Ising model with translation-invariant interactions.

Let us consider one more example. As before, let $M = \mathbf{T}$, and let $\mu_0$ be any probability on $\mathbf{T}$. Let $V : \mathbf{T} \times \mathbf{T} \to \mathbf{R}$ be a continuous bounded function that can be written as

$$V(x,y) = \sum_{k=1}^{\infty} a_k e_k(x) e_k(y), \qquad x, y \in \mathbf{T},$$

with $a_k < \frac{1}{2}$, $k \in \mathbf{N}$ and $\{e_k\}_{k\in\mathbf{N}}$ an ONB of $L_0^2(d\mu_0)$. [This implies in particular that $\int V(x,y)\mu_0(dy) = 0$ for any $x \in \mathbf{T}$.] Let $\Phi(R) = \iint V(x,y)R(dx)R(dy)$ for positive measures $R$ on $\mathbf{T}$. (The corresponding entropy function is again the relative entropy with respect to $\mu_0$.) It is easy to check that $\mu_0$ maximizes $\Phi - h$. Therefore, the space $H$ is given by $H = L_0^2(d\mu_0)$. So the eigenvalues and the corresponding eigenfunctions of $D^2\Phi(\mu_0)|_{H\times H}$ are $\{2a_k\}_{k\in\mathbf{N}}$ and $\{e_k\}_{k\in\mathbf{N}}$. Therefore, our assumptions are satisfied if, in addition, there exists a $p > 1$ such that

$$\sum_{k\in\mathbf{N} : a_k > 0} a_k^{1/p} e_k^{\otimes 2} : \mathbf{T} \times \mathbf{T} \to \mathbf{R} \text{ is bounded.}$$

The proof of Theorem 1.3 is given in Sections 2–5.

**2. Preparation.** Let us set $\lambda \equiv \Phi(x^*) - h(x^*) = \sup_{x\in B}\{\Phi(x) - h(x)\}$, for simplicity. Let $R_5(x,\cdot)$ denote the fifth remainder of the Taylor expansion of $\Phi(x + \cdot)$ at $x$, that is,

$$R_5(x,y) = \Phi(x+y) - \Phi(x) - \sum_{i=1}^{4} D^i\Phi(x)(y,\ldots,y) \qquad \text{for any } x, y \in B.$$

Then we have by [8] or [21] that for any $\varepsilon > 0$,

$$\begin{aligned}
&n(e^{-\lambda n}Z_n - \det(I_H - \Psi_2)^{-1/2}) \\
&\quad = n\Bigg(E^{\nu_0^{\otimes\infty}}\Bigg[\exp\Bigg(\frac{n}{2}\Psi_2\bigg(\frac{S_n}{n},\frac{S_n}{n}\bigg) + \frac{n}{3!}\Psi_3\bigg(\frac{S_n}{n},\frac{S_n}{n},\frac{S_n}{n}\bigg) \\
&\qquad\qquad + \frac{n}{4!}\Psi_4\bigg(\frac{S_n}{n},\frac{S_n}{n},\frac{S_n}{n},\frac{S_n}{n}\bigg) + nR_5\bigg(x^*,\frac{S_n}{n}\bigg)\Bigg)\Bigg), \\
&\qquad\qquad \bigg\|\frac{S_n}{n}\bigg\|_B < \varepsilon\Bigg] \\
&\qquad\qquad - \det(I_H - \Psi_2)^{-1/2}\Bigg) \\
&\qquad + ne^{-\lambda n}E^{\mu^{\otimes\infty}}\bigg[\exp\bigg(n\Phi\bigg(\frac{S_n}{n}\bigg)\bigg), \bigg\|\frac{S_n}{n} - x^*\bigg\|_B > \varepsilon\bigg],
\end{aligned}$$

and the second term on the right-hand side converges to 0 exponentially as $n \to \infty$, by using the large deviation principle (see, e.g., [21]). So we only need to deal with the first term on the right-hand side. We can rewrite it as

$$n\bigg(E^{\nu_0^{\otimes\infty}}\bigg[\exp\bigg(\frac{n}{2}\Psi_2\bigg(\frac{S_n}{n},\frac{S_n}{n}\bigg)\bigg), \bigg\|\frac{S_n}{n}\bigg\|_B < \varepsilon\bigg] - \det(I_H - \Psi_2)^{-1/2}\bigg) \tag{2.1}$$

$$\begin{aligned}
&+ nE^{\nu_0^{\otimes\infty}}\bigg[\exp\bigg(\frac{n}{2}\Psi_2\bigg(\frac{S_n}{n},\frac{S_n}{n}\bigg)\bigg) \\
&\qquad \times \bigg(\exp\bigg(\frac{n}{3!}\Psi_3\bigg(\frac{S_n}{n},\frac{S_n}{n},\frac{S_n}{n}\bigg)\bigg) - 1\bigg), \bigg\|\frac{S_n}{n}\bigg\|_B < \varepsilon\bigg]
\end{aligned} \tag{2.2}$$

$$\begin{aligned}
&+ E^{\nu_0^{\otimes\infty}}\bigg[\exp\bigg(\frac{n}{2}\Psi_2\bigg(\frac{S_n}{n},\frac{S_n}{n}\bigg) + \frac{n}{3!}\Psi_3\bigg(\frac{S_n}{n},\frac{S_n}{n},\frac{S_n}{n}\bigg)\bigg) \\
&\qquad \times n\bigg(\exp\bigg(\frac{n}{4!}\Psi_4\bigg(\frac{S_n}{n},\frac{S_n}{n},\frac{S_n}{n},\frac{S_n}{n}\bigg) \\
&\qquad\qquad + nR_5\bigg(x^*,\frac{S_n}{n}\bigg)\bigg) - 1\bigg), \bigg\|\frac{S_n}{n}\bigg\|_B < \varepsilon\bigg].
\end{aligned} \tag{2.3}$$

We will work with (2.1) in Section 3, (2.2) in Section 5 and (2.3) in Section 4, respectively.

**3. Second order.** In this section, we are going to give the asymptotic expansion of the term (2.1) for $n \to \infty$. The result is in fact stronger than what is needed for the proof of our Theorem 1.3, and of interest in itself (cf. Theorem 3.14). The basic idea is to first use the fact that the Laplace transform of a Gaussian measure is exp(quadratic form), then to use the independence of $X_k$, $k \in \mathbf{N}$, to discuss the a.s.-convergence and the dominations, which then implies the $L^1$-convergence.

In general, let $(B,\|\cdot\|_B)$ be a Banach space and let $X_n$, $n\in\mathbf{N}$ be a sequence of i.i.d. $B$-valued random variables with mean 0. Let $\nu_0$ denote their common distribution, and we suppose that the following assumption is satisfied.

(H1) There exists a constant $K_1>0$ such that

$$\int_B \exp(K_1\|x\|_B^2)\nu_0(dx)<\infty.$$

We remark that (H1) is equivalent to our Assumption A1.

Let $\Psi_2 : B\times B\to\mathbf{R}$ be a bilinear, symmetric, bounded map satisfying Assumptions A4 and A6 in Section 1.

For some $\varepsilon>0$ small enough, we want to know the asymptotic expansion of

$$E^{\nu_0^{\otimes\infty}}\left[\exp\left(\frac{n}{2}\Psi_2\left(\frac{S_n}{n},\frac{S_n}{n}\right)\right),\left\|\frac{S_n}{n}\right\|_B<\varepsilon\right],\qquad n\to\infty. \tag{3.1}$$

Let $H$ be the Hilbert space with norm given by $\|\int_B u(y)y\nu_0(dy)\|_H^2=\int u^2\,d\nu_0$.

The following two propositions follow easily from [21].

Proposition 3.1. *For any $\Psi : B\times B\to\mathbf{R}$ which is bilinear, symmetric, continuous, and satisfies the condition that all of the eigenvalues of $\Psi|_{H\times H}$ are strictly smaller than 1, there exists an $\varepsilon_0>0$ such that for any $\varepsilon\in(0,\varepsilon_0]$,*

$$\sup_{n\in\mathbf{N}} E^{\nu_0^{\otimes\infty}}\left[\exp\left(\frac{n}{2}\Psi\left(\frac{S_n}{n},\frac{S_n}{n}\right)\right),\left\|\frac{S_n}{n}\right\|_B<\varepsilon\right]<\infty.$$

Proposition 3.2. *For any $\Psi : B\times B\to\mathbf{R}$ which is symmetric, bilinear, continuous, there exists a $\delta_1>0$ such that*

$$\sup_{n\in\mathbf{N}} E^{\nu_0^{\otimes\infty}}\left[\exp\left(\delta_1\Psi\left(\frac{S_n}{\sqrt{n}},\frac{S_n}{\sqrt{n}}\right)\right)\right]<\infty.$$

By our Assumption A4, $a_k<1$ for all $k\in\mathbf{N}$. Also, we have $a_k\to0$ as $k\to\infty$, since $\Psi_2|_{H\times H}$ is a Hilbert–Schmidt operator by [21] (and even nuclear under our Assumption A6). So $a_k, k\in\mathbf{N}$, are uniformly separated from 1. Therefore, there exists a $p_0>1$ such that $p_0\cdot a_k<1$ for any $k\in\mathbf{N}$. Let $q_0>1$ be such that $\frac{1}{p_0}+\frac{1}{q_0}=1$. Let $N_0\in\mathbf{N}$ be (large enough) so that $q_0\cdot\delta_{N_0}<\delta_1$, where $\delta_N$ is the sequence of positive numbers that converges to 0 which appeared in Assumption A6, and $\delta_1>0$ is the constant which appeared in Proposition 3.2 applied to $\Psi=\widetilde{\Psi}_2$.

For any $N\in\mathbf{N}$, define $\|\cdot\|_{H_N}$ by $\|x\|_{H_N}\equiv\sum_{k:\,k\le N,\,a_k>0}e_k(x)^2$, $x\in B$. Then we have the following.

LEMMA 3.3. *For any $N \geq N_0$, there exists a constant $\varepsilon_0 > 0$ such that for any $\varepsilon \in (0, \varepsilon_0]$,*

$$\sup_{n \in \mathbf{N}} E^{\nu_0^{\otimes\infty}}\left[\exp\left(\frac{n}{2}\Psi_2\left(\frac{S_n}{n}, \frac{S_n}{n}\right)\right), \left\|\frac{S_n}{n}\right\|_{H_N} < \varepsilon\right] < \infty.$$

*Also, for any $\varepsilon_1 > 0$,*

$$E^{\nu_0^{\otimes\infty}}\left[\exp\left(\frac{n}{2}\Psi_2\left(\frac{S_n}{n}, \frac{S_n}{n}\right)\right), \left\{\left\|\frac{S_n}{n}\right\|_B > \varepsilon_1\right\} \cap \left\{\left\|\frac{S_n}{n}\right\|_{H_N} < \varepsilon\right\}\right]$$

*converges to 0 exponentially as $n \to \infty$.*

PROOF. The first assertion is easily proven since

$$\begin{aligned} &E^{\nu_0^{\otimes\infty}}\left[\exp\left(\frac{n}{2}\Psi_2\left(\frac{S_n}{n}, \frac{S_n}{n}\right)\right), \left\|\frac{S_n}{n}\right\|_{H_N} < \varepsilon\right] \\ &\quad \leq E^{\nu_0^{\otimes\infty}}\left[\exp\left(\frac{p_0}{2}\sum_{k:\, k \leq N, a_k > 0} a_k e_k\left(\frac{S_n}{\sqrt{n}}\right)^2\right), \left\|\frac{S_n}{n}\right\|_{H_N} < \varepsilon\right]^{1/p_0} \\ &\qquad \times E^{\nu_0^{\otimes\infty}}\left[\exp\left(q_0 \delta_N \widetilde{\Psi}_2\left(\frac{S_n}{\sqrt{n}}, \frac{S_n}{\sqrt{n}}\right)\right)\right]^{1/q_0}, \end{aligned}$$

and the first term above is bounded for $n \in \mathbf{N}$ by Proposition 3.1 with $B$, $H$ and $\Psi$ replaced by $B_N$, $B_N$ and $p_0 \cdot \sum_{k:\, k \leq N, a_k > 0} a_k e_k^{\otimes 2}$, respectively, and the second term above is bounded for $n \in \mathbf{N}$ by Proposition 3.2 and the definition of $N_0$.

For the second assertion, choose $r > 1$ so that $r \cdot p_0 \cdot a_k < 1$ for all $k \in \mathbf{N}$, and let $s > 1$ be such that $\frac{1}{r} + \frac{1}{s} = 1$. Then

$$\begin{aligned} &E^{\nu_0^{\otimes\infty}}\left[\exp\left(\frac{n}{2}\Psi_2\left(\frac{S_n}{n}, \frac{S_n}{n}\right)\right), \left\{\left\|\frac{S_n}{n}\right\|_B > \varepsilon_1\right\} \cap \left\{\left\|\frac{S_n}{n}\right\|_{H_N} < \varepsilon\right\}\right] \\ &\quad \leq E^{\nu_0^{\otimes\infty}}\left[\exp\left(\frac{r}{2}\Psi_2\left(\frac{S_n}{\sqrt{n}}, \frac{S_n}{\sqrt{n}}\right)\right), \left\|\frac{S_n}{n}\right\|_{H_N} < \varepsilon\right]^{1/r} P^{\nu_0^{\otimes\infty}}\left(\left\|\frac{S_n}{n}\right\|_B > \varepsilon_1\right)^{1/s}. \end{aligned}$$

The first term is bounded for $n \in \mathbf{N}$ by our first assertion applied to $r \cdot \Psi_2$, and the second term decays exponentially as $n \to \infty$, by the large deviation principle, the properties of the entropy function $h$ and the fact that $\nu_0$ is 0-centered. This completes the proof of our lemma. □

We have by Lemma 3.3 that there exists an $N_0 \in \mathbf{N}$ such that for any $N \geq N_0$, there exists an $\varepsilon_0 > 0$ such that for any $\varepsilon \in (0, \varepsilon_0]$, the asymptotic expansion of (3.1) is the same as the asymptotic expansion of

$$E^{\nu_0^{\otimes\infty}}\left[\exp\left(\frac{n}{2}\Psi_2\left(\frac{S_n}{n}, \frac{S_n}{n}\right)\right), \left\|\frac{S_n}{n}\right\|_{H_N} < \varepsilon\right], \tag{3.2}$$

for $n \to \infty$.

From now on, we let $M \geq N_0$ be chosen and fixed. Let $H_1$ be the Hilbert space given by

$$H_1 = \left\{ y = \sum_{k=1}^{\infty} e_k(y) f_k ; \|y\|_{H_1}^2 = \sum_{k=1}^{\infty} |a_k| e_k(y)^2 < \infty \right\}.$$

Since $\sum_{k=1}^{\infty} |a_k| < \infty$ by our Assumption A6, we have that there exists an $H_1$-valued Gaussian random variable $Y$ such that $u(Y) \sim N(0, \|u\|_{H^*}^2)$ for any $u \in H^*$ ("$\sim$" meaning equality in law). Let $\hat{Y^M} \equiv \sum_{k\,:\,k \leq M,\, a_k > 0} e_k(Y) f_k$.

It is easy to see, using the Fourier transform of the Gaussian measure, that for any $x \in B$,

$$\exp\left( \tfrac{1}{2} \sum_{k=1}^{\infty} a_k e_k(x)^2 \right) = E^Y \left[ \exp\left( \sum_{k=1}^{\infty} \sqrt{a_k} e_k(x) e_k(Y) \right) \right]. \tag{3.3}$$

We are going to use this fact to give the asymptotic expansion of (3.2) as $n \to \infty$.

First, by Assumption A6, (3.3) and Fubini's theorem, we have that

$$\begin{aligned} & E^{\nu_0^{\otimes\infty}} \left[ \exp\left( \frac{n}{2} \Psi_2\left( \frac{S_n}{n}, \frac{S_n}{n} \right) \right), \left\| \frac{S_n}{n} \right\|_{H_M} < \varepsilon \right] \\ &= E^{\nu_0^{\otimes\infty}} \left[ \exp\left( \frac{n}{2} \sum_{k=1}^{\infty} a_k e_k \left( \frac{S_n}{n} \right)^2 \right), \left\| \frac{S_n}{n} \right\|_{H_M} < \varepsilon \right] \\ &= E^{\nu_0^{\otimes\infty}} \left[ E^Y \left[ \exp\left( \sum_{k=1}^{\infty} \sqrt{a_k} e_k \left( \frac{S_n}{\sqrt{n}} \right) e_k(Y) \right) \right], \left\| \frac{S_n}{n} \right\|_{H_M} < \varepsilon \right] \\ &= E^Y \left[ E^{\nu_0^{\otimes\infty}} \left[ \exp\left( \sum_{k=1}^{\infty} \sqrt{a_k} e_k \left( \frac{S_n}{\sqrt{n}} \right) e_k(Y) \right), \left\| \frac{S_n}{n} \right\|_{H_M} < \varepsilon \right] \right]. \end{aligned} \tag{3.4}$$

Proposition 3.4. *There exists an $\varepsilon_0 > 0$ such that for any $\varepsilon \in (0, \varepsilon_0]$ and any $\delta > 0$,*

$$E^Y \left[ E^{\nu_0^{\otimes\infty}} \left[ \exp\left( \sum_{k=1}^{\infty} \sqrt{a_k} e_k \left( \frac{S_n}{\sqrt{n}} \right) e_k(Y) \right), \left\| \frac{S_n}{n} \right\|_{H_M} < \varepsilon \right], \|\hat{Y^M}\|_{H_1} > \sqrt{n}\varepsilon\delta \right]$$

*converges to 0 exponentially as $n \to \infty$.*

PROOF. We have by Hölder's inequality that for any $p, q > 1$ such that $\frac{1}{p} + \frac{1}{q} = 1$,

$$
\begin{aligned}
&\Bigg| E^Y \Bigg[ E^{\nu_0^{\otimes\infty}} \Bigg[ \exp\Bigg( \sum_{k=1}^{\infty} \sqrt{a_k} e_k\Big(\frac{S_n}{\sqrt{n}}\Big) e_k(Y) \Bigg), \Big\| \frac{S_n}{n} \Big\|_{H_M} < \varepsilon \Bigg], \\
&\qquad\qquad\qquad\qquad \|Y^{\hat{M}}\|_{H_1} > \sqrt{n}\varepsilon\delta \Bigg] \Bigg| \\
&\quad = E^{\nu_0^{\otimes\infty}} \Bigg[ \exp\Bigg( \frac{1}{2} \sum_{k\,:\, k>M \text{ or } a_k<0} a_k e_k\Big(\frac{S_n}{\sqrt{n}}\Big)^2 \Bigg) \\
&\qquad\quad \times E^Y \Bigg[ \exp\Bigg( \sum_{k\,:\, k\le M \text{ and } a_k>0} \sqrt{a_k} e_k\Big(\frac{S_n}{\sqrt{n}}\Big) e_k(Y) \Bigg), \\
&\qquad\qquad\qquad\qquad \|Y^{\hat{M}}\|_{H_1} > \sqrt{n}\varepsilon\delta \Bigg], \Big\| \frac{S_n}{n} \Big\|_{H_M} < \varepsilon \Bigg] \\
&\quad \le E^{\nu_0^{\otimes\infty}} \Bigg[ \exp\Bigg( \frac{1}{2} \sum_{k\,:\, k>M \text{ or } a_k<0} a_k e_k\Big(\frac{S_n}{\sqrt{n}}\Big)^2 \Bigg) \\
&\qquad\quad \times E^Y \Bigg[ \exp\Bigg( p \cdot \sum_{k\,:\, k\le M \text{ and } a_k>0} \sqrt{a_k} e_k\Big(\frac{S_n}{\sqrt{n}}\Big) e_k(Y) \Bigg) \Bigg]^{1/p} \\
&\qquad\qquad\qquad \times P(\|Y^{\hat{M}}\|_{H_1} > \sqrt{n}\varepsilon\delta)^{1/q}, \Big\| \frac{S_n}{n} \Big\|_{H_M} < \varepsilon \Bigg] \\
&\quad = E^{\nu_0^{\otimes\infty}} \Bigg[ \exp\Bigg( \frac{1}{2} \sum_{k=1}^{\infty} a_k e_k\Big(\frac{S_n}{\sqrt{n}}\Big)^2 \\
&\qquad\qquad + \frac{1}{2}(p-1) \sum_{k\,:\, k\le M, a_k>0} a_k e_k\Big(\frac{S_n}{\sqrt{n}}\Big)^2 \Bigg), \Big\| \frac{S_n}{n} \Big\|_{H_M} < \varepsilon \Bigg] \\
&\qquad \times P(\|Y^{\hat{M}}\|_{H_1} > \sqrt{n}\varepsilon\delta)^{1/q}.
\end{aligned}
$$

The first factor in the latter expression is bounded for $n \in \mathbf{N}$ if $p > 1$ and $\varepsilon > 0$ are small enough, by Assumption A6 and Lemma 3.3, and the second factor decays exponentially as $n \to \infty$ for any $\delta > 0$. This gives our assertion. □

Proposition 3.5. *There exists a* $\delta_0 > 0$ *such that for any* $\delta \in (0, \delta_0]$ *and any* $\varepsilon > 0$,

$$E^Y\left[E^{\nu_0^{\otimes\infty}}\left[\exp\left(\sum_{k=1}^{\infty}\sqrt{a_k}e_k\left(\frac{S_n}{\sqrt{n}}\right)e_k(Y)\right), \left\|\frac{S_n}{n}\right\|_{H_M} > \varepsilon\right], \|Y^{\hat{M}}\|_{H_1} < \sqrt{n}\varepsilon\delta\right]$$

*converges to* 0 *exponentially as* $n \to \infty$.

Proof. First notice that

$$\left|\sum_{k:\, k\le M,\, a_k>0}\sqrt{a_k}e_k\left(\frac{S_n}{\sqrt{n}}\right)e_k(Y)\right| \le \left\|\frac{S_n}{\sqrt{n}}\right\|_{H_M}\cdot\|Y^{\hat{M}}\|_{H_1} \le \delta\left\|\frac{S_n}{\sqrt{n}}\right\|_{H_M}^2$$
$$= \delta\sum_{k:\, k\le M,\, a_k>0} e_k\left(\frac{S_n}{\sqrt{n}}\right)^2$$

on the set $\{\|\frac{S_n}{n}\|_{H_M} > \varepsilon, \|Y^{\hat{M}}\|_{H_1} < \sqrt{n}\varepsilon\delta\}$, and

$$\sum_{k>M,\, a_k>0} a_k e_k^{\otimes 2} \le \delta_M\widetilde{\Psi}_2$$

by Assumption A6. So we have that

$$\left|E^Y\left[E^{\nu_0^{\otimes\infty}}\left[\exp\left(\sum_{k=1}^{\infty}\sqrt{a_k}e_k\left(\frac{S_n}{\sqrt{n}}\right)e_k(Y)\right), \left\|\frac{S_n}{n}\right\|_{H_M} > \varepsilon\right],\right.\right.$$
$$\left.\left.\|Y^{\hat{M}}\|_{H_1} < \sqrt{n}\varepsilon\delta\right]\right|$$
$$= E^{\nu_0^{\otimes\infty}}\left[\exp\left(\frac{1}{2}\sum_{k:\, k>M \text{ or } a_k<0} a_k e_k\left(\frac{S_n}{\sqrt{n}}\right)^2\right)\right.$$
$$\times E^Y\left[\exp\left(\sum_{k:\, k\le M \text{ and } a_k>0}\sqrt{a_k}e_k\left(\frac{S_n}{\sqrt{n}}\right)e_k(Y)\right),\right.$$
$$\left.\left.\|Y^{\hat{M}}\|_{H_1} < \sqrt{n}\varepsilon\delta\right], \left\|\frac{S_n}{n}\right\|_{H_M} > \varepsilon\right]$$
$$\le E^{\nu_0^{\otimes\infty}}\left[\exp\left(\delta\sum_{k:\, k\le M,\, a_k>0} a_k e_k\left(\frac{S_n}{\sqrt{n}}\right)^2 + \frac{1}{2}\delta_M\widetilde{\Psi}_2\left(\frac{S_n}{\sqrt{n}}, \frac{S_n}{\sqrt{n}}\right)\right),\right.$$
$$\left.\left\|\frac{S_n}{n}\right\|_{H_M} > \varepsilon\right]$$
$$\le E^{\nu_0^{\otimes\infty}}\left[\exp\left(2\delta\sum_{k:\, k\le M,\, a_k>0} a_k e_k\left(\frac{S_n}{\sqrt{n}}\right)^2 + \delta_M\widetilde{\Psi}_2\left(\frac{S_n}{\sqrt{n}}, \frac{S_n}{\sqrt{n}}\right)\right)\right]^{1/2}$$

$$\times P\bigg(\bigg\|\frac{S_n}{n}\bigg\|_{H_M} > \varepsilon\bigg)^{1/2}.$$

Now, our assertion follows easily by Proposition 3.2 and the fact that $P(\|\frac{S_n}{n}\|_{H_M} > \varepsilon) \to 0$ exponentially as $n \to \infty$. □

In particular, we have the following.

Proposition 3.6. *There exist constants $M \in \mathbf{N}$, $\varepsilon_0 > 0$ and $\delta_0 > 0$ such that for any $\varepsilon \in (0, \varepsilon_0]$ and any $\delta \in (0, \delta_0]$,*

$$\sup_{n \in \mathbf{N}} E^{\nu_0^{\otimes\infty}}\bigg[E^Y\bigg[\exp\bigg(\sum_{k=1}^{\infty}\sqrt{a_k}e_k\bigg(\frac{S_n}{\sqrt{n}}\bigg)e_k(Y)\bigg), \|Y^{\hat{M}}\|_{H_1} < \sqrt{n}\varepsilon\delta\bigg]\bigg] < \infty.$$

In the same way, we get the following.

Proposition 3.7. *There exist constants $p > 1$, $M \in \mathbf{N}$, $\varepsilon_0 > 0$ and $\delta_0 > 0$ such that for any $\varepsilon \in (0, \varepsilon_0]$ and any $\delta \in (0, \delta_0]$,*

$$\sup_{n \in \mathbf{N}} E^{\nu_0^{\otimes\infty}}\bigg[E^Y\bigg[\bigg|\exp\bigg(\sum_{k=1}^{\infty}\sqrt{a_k}e_k\bigg(\frac{S_n}{\sqrt{n}}\bigg)e_k(Y)\bigg)\bigg|^p, \|Y^{\hat{M}}\|_{H_1} < \sqrt{n}\varepsilon\delta\bigg]\bigg] < \infty.$$

By (3.4) and Propositions 3.4, 3.5 and 3.7, we have that there exist constants $\varepsilon_0 > 0$ and $\delta_0 > 0$ such that for any $\varepsilon \in (0, \varepsilon_0]$ and any $\delta \in (0, \delta_0]$, the asymptotic expansion of (3.2) is the same as the asymptotic expansion of

$$\begin{aligned} &E^Y\bigg[E^{\nu_0^{\otimes\infty}}\bigg[\exp\bigg(\sum_{k=1}^{\infty}\sqrt{a_k}e_k\bigg(\frac{S_n}{\sqrt{n}}\bigg)e_k(Y)\bigg)\bigg], \|Y^{\hat{M}}\|_{H_1} < \sqrt{n}\varepsilon\delta\bigg] \\ &\quad = E^Y\bigg[E^{\nu_0}\bigg[\exp\bigg(\frac{1}{\sqrt{n}}\sum_{k=1}^{\infty}\sqrt{a_k}e_k(X_1)e_k(Y)\bigg)\bigg]^n, \|Y^{\hat{M}}\|_{H_1} < \sqrt{n}\varepsilon\delta\bigg], \end{aligned} \tag{3.5}$$

for $n \to \infty$.

For any $n \in \mathbf{N}$ and any $\xi > 0$, let $B_{n,\xi}$ be the set given by

$$B_{n,\xi} := \bigg\{Y : \bigg|E^{\nu_0}\bigg[\exp\bigg(\frac{1}{\sqrt{n}}\sum_{k=1}^{\infty}\sqrt{a_k}e_k(X_1)e_k(Y)\bigg)\bigg] - 1\bigg| < \xi\bigg\}. \tag{3.6}$$

Note that

$$\bigg|E^{\nu_0}\bigg[\exp\bigg(\frac{1}{\sqrt{n}}\sum_{k=1}^{\infty}\sqrt{a_k}e_k(X_1)e_k(Y)\bigg)\bigg] - 1\bigg|$$

$$\leq \frac{1}{n} E^{\nu_0}\left[\exp\left(\left|\frac{1}{\sqrt{n}}\sum_{k=1}^{\infty}\sqrt{a_k}e_k(X_1)e_k(Y)\right|\right)\right.$$
$$\left.\times\left|\sum_{k=1}^{\infty}\sqrt{a_k}e_k(X_1)e_k(Y)\right|^2\right],$$

and for any $q>1$, there exists a $N(q)\in\mathbf{N}$ such that

$$\sup_{n>N(q)} E^Y\left[\left(E^{\nu_0}\left[\exp\left(\left|\frac{1}{\sqrt{n}}\sum_{k=1}^{\infty}\sqrt{a_k}e_k(X_1)e_k(Y)\right|\right)\right.\right.\right.$$
$$\left.\left.\left.\times\left|\sum_{k=1}^{\infty}\sqrt{a_k}e_k(X_1)e_k(Y)\right|^2\right]\right)^q\right]<\infty.$$

Therefore,

$$P(B_{n,\xi}^c) \leq P^Y\left(E^{\nu_0}\left[\exp\left(\left|\frac{1}{\sqrt{n}}\sum_{k=1}^{\infty}\sqrt{a_k}e_k(X_1)e_k(Y)\right|\right)\right.\right.$$
$$\left.\left.\times\left|\sum_{k=1}^{\infty}\sqrt{a_k}e_k(X_1)e_k(Y)\right|^2\right]>\xi n\right)$$
$$\leq (\xi n)^{-q}E^Y\left[\left(E^{\nu_0}\left[\exp\left(\left|\frac{1}{\sqrt{n}}\sum_{k=1}^{\infty}\sqrt{a_k}e_k(X_1)e_k(Y)\right|\right)\right.\right.\right.$$
$$\left.\left.\left.\times\left|\sum_{k=1}^{\infty}\sqrt{a_k}e_k(X_1)e_k(Y)\right|^2\right]\right)^q\right].$$

So we have the following.

PROPOSITION 3.8. *$P(B_{n,\xi}^c)\to 0$ as $n\to\infty$ faster than any polynomial order. Therefore,*

$$E^Y\left[E^{\nu_0}\left[\exp\left(\frac{1}{\sqrt{n}}\sum_{k=1}^{\infty}\sqrt{a_k}e_k(X_1)e_k(Y)\right)\right]^n,\{\|\hat{Y^M}\|_{H_1}<\sqrt{n}\varepsilon\delta\}\cap B_{n,\xi}^c\right]\to 0$$

*faster than any polynomial order.*

PROOF. The first assertion is already proven. The second one follows then, using Hölder's inequality and Proposition 3.7. □

By Proposition 3.8, we have that for any $\xi>0$, the asymptotic expansion of (3.5) for $n\to\infty$ is the same as the corresponding asymptotic expansion

of

$$E^Y\left[E^{\nu_0}\left[\exp\left(\frac{1}{\sqrt{n}}\sum_{k=1}^{\infty}\sqrt{a_k}e_k(X_1)e_k(Y)\right)\right]^n,\right.\tag{3.7}$$
$$\left.\{\|Y^{\hat{M}}\|_{H_1}<\sqrt{n}\varepsilon\delta\}\cap B_{n,\xi}\right].$$

We can take, for example, $\xi=\frac{1}{2}$ in the definition (3.6) of $B_{n,\xi}$. Let

$$Z\equiv E^{\nu_0}\left[\exp\left(\frac{1}{\sqrt{n}}\sum_{k=1}^{\infty}\sqrt{a_k}e_k(X_1)e_k(Y)\right)\right]-1.$$

Then we have the following.

Proposition 3.9. *For any $N\geq 2$,*

$$n^{N/2}\left|Z-\sum_{j=2}^{N}\frac{1}{j!}n^{-j/2}E^{\nu_0}\left[\left(\sum_{k=1}^{\infty}\sqrt{a_k}e_k(X_1)e_k(Y)\right)^j\right]\right|\to 0\qquad \textit{as } n\to\infty$$

*for a.e.-$Y$, and for any $q>1$, there exists an $n_0\in\mathbf{N}$ such that for any $n\geq n_0$, the left-hand side above is dominated by an $L^q$ (with respect to the distribution of the $Y$)-random variable.*

Proof. We notice that $|e^x-\sum_{j=0}^{N}\frac{x^j}{j!}|\leq e^{|x|}\frac{|x|^{N+1}}{(N+1)!}$ for any $x\in\mathbf{C}$.

First, as we claimed before, $|\sum_{k=1}^{\infty}\sqrt{a_k}e_k(X_1)e_k(Y)|<\infty$, so we have that

$$n^{N/2}\left(\exp\left(\frac{1}{\sqrt{n}}\sum_{k=1}^{\infty}\sqrt{a_k}e_k(X_1)e_k(Y)\right)\right.$$
$$\left.-1-\sum_{j=1}^{N}\frac{1}{j!}n^{-j/2}\left(\sum_{k=1}^{\infty}\sqrt{a_k}e_k(X_1)e_k(Y)\right)^j\right)\to 0$$

as $n\to\infty$ for a.e.-$(X_1,Y)$. Also, for any $\eta>0$, there exists an $n_0\in\mathbf{N}$ (actually, we can take, e.g., $n_0=[\frac{1}{\eta^2}]+1$) such that for any $n\geq n_0$,

$$\left|n^{N/2}\left(\exp\left(\frac{1}{\sqrt{n}}\sum_{k=1}^{\infty}\sqrt{a_k}e_k(X_1)e_k(Y)\right)\right.\right.$$
$$\left.\left.-1-\sum_{j=1}^{N}\frac{1}{j!}n^{-j/2}\left(\sum_{k=1}^{\infty}\sqrt{a_k}e_k(X_1)e_k(Y)\right)^j\right)\right|$$
$$\leq\frac{1}{\sqrt{n}}\exp\left(\left|\frac{1}{\sqrt{n}}\sum_{k=1}^{\infty}\sqrt{a_k}e_k(X_1)e_k(Y)\right|\right)\frac{1}{(N+1)!}\left|\sum_{k=1}^{\infty}\sqrt{a_k}e_k(X_1)e_k(Y)\right|^{N+1}$$

$$\leq \frac{1}{\sqrt{n}} \frac{1}{\eta^{N+1}} \exp\left(2\eta\left|\sum_{k=1}^{\infty} \sqrt{a_k} e_k(X_1) e_k(Y)\right|\right).$$

Also, since

$$\begin{aligned}
&\exp\left(2\eta\left|\sum_{k=1}^{\infty} \sqrt{a_k} e_k(X_1) e_k(Y)\right|\right) \\
&\quad \leq \left(\exp\left(2\eta \sum_{k\,:\,a_k>0} \sqrt{a_k} e_k(X_1) e_k(Y)\right)\right. \\
&\qquad \left. + \exp\left(-2\eta \sum_{k\,:\,a_k>0} \sqrt{a_k} e_k(X_1) e_k(Y)\right)\right) \\
&\qquad \times \left(\exp\left(2\eta \sum_{k\,:\,a_k<0} \sqrt{|a_k|} e_k(X_1) e_k(Y)\right)\right. \\
&\qquad \left. + \exp\left(-2\eta \sum_{k\,:\,a_k<0} \sqrt{|a_k|} e_k(X_1) e_k(Y)\right)\right),
\end{aligned}$$

we have that

$$\begin{aligned}
&E^Y\left[E^{\nu_0}\left[\exp\left(2\eta\left|\sum_{k=1}^{\infty} \sqrt{a_k} e_k(X_1) e_k(Y)\right|\right)\right]^q\right] \\
&\quad \leq E^{\nu_0}\left[2\exp\left(2q^2\eta^2 \sum_{k\,:\,a_k>0} a_k e_k(X_1)^2\right) 2\exp\left(2q^2\eta^2 \sum_{k\,:\,a_k<0} a_k e_k(X_1)^2\right)\right] \\
&\quad = 4E^{\nu_0}\left[\exp\left(2q^2\eta^2 \sum_{k=1}^{\infty} |a_k| e_k(X_1)^2\right)\right] \\
&\quad \leq 4E^{\nu_0}\left[\exp\left(2q^2\eta^2 \left(\max_{k\in\mathbf{N}} |a_k|\right) \|X_1\|_B^2\right)\right] < \infty
\end{aligned}$$

if $\eta > 0$ is small enough [so that $2q^2\eta^2(\max_{k\in\mathbf{N}} |a_k|) \leq K_1$].

Therefore, we get our assertion by the dominated convergence theorem. □

When $N = 2$, since

$$E^{\nu_0}\left[\left(\sum_{k=1}^{\infty} \sqrt{a_k} e_k(X_1) e_k(Y)\right)^2\right] = \sum_{k=1}^{\infty} a_k e_k(Y)^2,$$

Proposition 3.9 gives us the following.

COROLLARY 3.10.

$$n^{1/2}\left|n\log(1+Z)-\tfrac{1}{2}\sum_{k=1}^{\infty}a_k e_k(Y)^2\right| \tag{3.8}$$

*is bounded in $L^q$ (with respect to the distribution of $Y$) for any $q>1$ on the set $B_{n,1/2}$.*

PROOF. We notice that there exists a constant $C>0$ such that

$$|\log(1+Z)-Z|\le C|Z|^2 \qquad \text{on } B_{n,1/2}=\{|Z|<\tfrac{1}{2}\}.$$

So

$$\begin{aligned}
&\text{(lhs) of (3.8)}\\
&\quad\le n^{1/2}|n\log(1+Z)-nZ|+n^{1/2}\left|nZ-\tfrac{1}{2}E^{\nu_0}\left[\left(\sum_{k=1}^{\infty}\sqrt{a_k}e_k(X_1)e_k(Y)\right)^2\right]\right|\\
&\quad\le Cn^{-1/2}|nZ|^2+n^{1/2}\left|nZ-\tfrac{1}{2}E^{\nu_0}\left[\left(\sum_{k=1}^{\infty}\sqrt{a_k}e_k(X_1)e_k(Y)\right)^2\right]\right|,
\end{aligned}$$

which is bounded in $L^q$ (with respect to the distribution of $Y$) for any $q>1$ by Proposition 3.9. □

PROPOSITION 3.11. *There exist constants $p_1>1$, $\varepsilon_0>0$ and $\delta_0>0$ such that for any $p<p_1$, any $\varepsilon\in(0,\varepsilon_0]$ any $\delta\in(0,\delta_0]$ and any $N\in\mathbf{N}$,*

$$\begin{aligned}
&\left|n^N\exp\left(\frac{1}{2}\sum_{k=1}^{\infty}a_k e_k(Y)^2\right)\right.\\
&\quad\times\left(\exp\left(n\log(1+Z)-\frac{1}{2}\sum_{k=1}^{\infty}a_k e_k(Y)^2\right)-1\right.\\
&\quad\left.\left.-\sum_{\ell=1}^{2N}\frac{1}{\ell!}\left(n\log(1+Z)-\frac{1}{2}\sum_{k=1}^{\infty}a_k e_k(Y)^2\right)^{\ell}\right)\right|\to 0
\end{aligned} \tag{3.9}$$

*for a.e.-$Y$, and the left-hand side above is $L^p$ (with respect to the distribution of $Y$)-bounded on the set $\{\|\hat{Y^M}\|_{H_1}<\sqrt{n}\varepsilon\delta\}\cap B_{n,1/2}$.*

PROOF. We notice that for any $N\in\mathbf{N}$, there exists a constant $C_N>0$ such that

$$\left|e^z-1-\sum_{\ell=1}^{2N}\frac{z^\ell}{\ell!}\right|\le C_N(|e^z|\vee 1)|z|^{2N+1} \qquad \text{for any } z\in\mathbf{C}.$$

Therefore,

$$\text{(lhs) of (3.9)} \le n^N \left( |e^{n\log(1+Z)}| + \left| \exp\left( \tfrac{1}{2} \sum_{k=1}^{\infty} a_k e_k(Y)^2 \right) \right| \right)$$
$$\times \left| n\log(1+Z) - \tfrac{1}{2} \sum_{k=1}^{\infty} a_k e_k(Y)^2 \right|^{2N+1}.$$

There exist constants $p > 1$ and $r > 1$ such that

$$\left( |e^{n\log(1+Z)}| + \left| \exp\left( \tfrac{1}{2} \sum_{k=1}^{\infty} a_k e_k(Y)^2 \right) \right| \right)$$

is bounded in $L^{p\cdot r}$ (with respect to the distribution of $Y$) on the set $\{\|\hat{Y^M}\|_{H_1} < \sqrt{n}\varepsilon\delta\}$. (The assertion for the first term is by Proposition 3.7, and the assertion for the second term is easy by Assumption A4.) Let $s > 1$ be such that $\frac{1}{r} + \frac{1}{s} = 1$. By Corollary 3.10,

$$n^{1/2} \left| n\log(1+Z) - \tfrac{1}{2} \sum_{k=1}^{\infty} a_k e_k(Y)^2 \right|$$

is bounded in $L^q$ (with respect to the distribution of $Y$) for any $q > 1$ on the set $B_{n,1/2}$; in particular, it is bounded in $L^{p\cdot s}$ (with respect to the distribution of $Y$) on the set $B_{n,1/2}$.

These give us our assertion. □

Proposition 3.12. *For any $N \in \mathbf{N}$,*

$$n^N \left| \sum_{\ell=1}^{2N} \frac{1}{\ell!} \left( n\log(1+Z) - \frac{1}{2} \sum_{k=1}^{\infty} a_k e_k(Y)^2 \right)^{\ell} \right. \tag{3.10}$$

$$\left. - \sum_{\ell=1}^{2N} \frac{1}{\ell!} \left( n \sum_{j=1}^{N+1} \frac{(-1)^{j-1}}{j} Z^j - \frac{1}{2} \sum_{k=1}^{\infty} a_k e_k(Y)^2 \right)^{\ell} \right| \to 0 \tag{3.11}$$

*for a.e.-$Y$, and the left-hand side above is $L^q$ (with respect to the distribution of $Y$)-bounded on the set $B_{n,1/2}$ for any $q > 1$.*

Proof. We notice that $|a^\ell - b^\ell| \le |a-b| \cdot \ell(|a|^{\ell-1} + |b|^{\ell-1})$ for any $a, b \in \mathbf{C}$ and any $\ell \in \mathbf{N}$. So our assertion is easy since by Proposition 3.9 and Corollary 3.10, $n\log(1+Z) - \frac{1}{2}\sum_{k=1}^{\infty} a_k e_k(Y)^2$ and $n\sum_{j=1}^{N+1} \frac{(-1)^j}{j} Z^j - \frac{1}{2}\sum_{k=1}^{\infty} a_k e_k(Y)^2$ are $L^q$ (with respect to the distribution of $Y$)-bounded on the set $B_{n,1/2}$ for any $q > 1$, and for any $N \in \mathbf{N}$, there exists a constant $C_N > 0$ such that

$$\left| n^N \left( n\log(1+Z) - n \sum_{j=1}^{N+1} \frac{(-1)^{j-1}}{j} Z^j \right) \right| \le C_N \frac{1}{n} |nZ|^{N+2}$$

on the set $B_{n,1/2}$, and $|nZ|$ is $L^q$ (with respect to the distribution of $Y$)-bounded by Proposition 3.9 for any $q > 1$. □

In the same way, by Proposition 3.9, we get the following.

PROPOSITION 3.13. *For any $N \in \mathbf{N}$,*

$$n^N \Bigg| \sum_{\ell=1}^{2N} \frac{1}{\ell!} \left( n \sum_{j=1}^{N+1} \frac{(-1)^{j-1}}{j} Z^j - \frac{1}{2} \sum_{k=1}^{\infty} a_k e_k(Y)^2 \right)^{\ell} - \sum_{\ell=1}^{2N} \frac{1}{\ell!} \left( \sum_{m=3}^{2N+2} n^{-(m/2)+1} \frac{1}{m!} E^{\nu_0} \left[ \left( \sum_{k=1}^{\infty} \sqrt{a_k} e_k(X_1) e_k(Y) \right)^m \right] + n \sum_{j=2}^{N+1} \frac{(-1)^{j-1}}{j} \times \left( \sum_{m=2}^{N+1} n^{-m/2} \frac{1}{m!} E^{\nu_0} \left[ \left( \sum_{k=1}^{\infty} \sqrt{a_k} e_k(X_1) e_k(Y) \right)^m \right] \right)^j \right)^{\ell} \Bigg| \to 0$$

*as $n \to \infty$ for a.e.-$Y$, and the left-hand side above is $L^q$ (with respect to the distribution of $Y$)-bounded for any $q > 1$.*

In conclusion, we have proven the following.

THEOREM 3.14. *There exist constants $\varepsilon_0 > 0$ and $\delta_0 > 0$ such that for any $\varepsilon \in (0, \varepsilon_0]$, any $\delta \in (0, \delta_0]$ and any $N \in \mathbf{N}$,*

$$n^N \Bigg\{ E^{\nu_0} \left[ \exp\left( \frac{1}{\sqrt{n}} \sum_{k=1}^{\infty} \sqrt{a_k} e_k(X_1) e_k(Y) \right) \right]^n - \exp\left( \frac{1}{2} \sum_{k=1}^{\infty} a_k e_k(Y)^2 \right) - \exp\left( \frac{1}{2} \sum_{k=1}^{\infty} a_k e_k(Y)^2 \right) \times \sum_{\ell=1}^{2N} \frac{1}{\ell!} \left( \sum_{m=3}^{2N+2} n^{-(m/2)+1} \frac{1}{m!} E^{\nu_0} \left[ \left( \sum_{k=1}^{\infty} \sqrt{a_k} e_k(X_1) e_k(Y) \right)^m \right] + n \sum_{j=2}^{N+1} \frac{(-1)^{j-1}}{j} \times \left( \sum_{m=2}^{N+1} n^{-m/2} \frac{1}{m!} E^{\nu_0} \left[ \left( \sum_{k=1}^{\infty} \sqrt{a_k} e_k(X_1) e_k(Y) \right)^m \right] \right)^j \right)^{\ell} \Bigg\} \to 0 \qquad \text{as } n \to \infty$$

*for a.e.-$Y$. Moreover, there exists a $p>1$ such that the left-hand side is $L^p$ (with respect to the distribution of $Y$)-bounded on the set $\{\|\hat{Y^M}\|_{H_1}\le \sqrt{n}\varepsilon\delta\}\cap B_{n,1/2}$ with respect to $n\in\mathbf{N}$. Therefore, we get that for any $\varepsilon\in(0,\varepsilon_0]$ and any $\delta\in(0,\delta_0]$,*

$$
\begin{aligned}
n^N\Biggl(E^{\nu_0^{\otimes\infty}}\Biggl[\exp\biggl(\frac{n}{2}\Psi_2\biggl(\frac{S_n}{n},\frac{S_n}{n}\biggr)\biggr),\biggl\|\frac{S_n}{n}\biggr\|_B<\varepsilon\Biggr]-E^Y\Biggl[\exp\Biggl(\frac{1}{2}\sum_{k=1}^{\infty}a_ke_k(Y)^2\Biggr)\Biggr]\\
-E^Y\Biggl[\exp\Biggl(\frac{1}{2}\sum_{k=1}^{\infty}a_ke_k(Y)^2\Biggr)\\
\times\sum_{\ell=1}^{2N}\frac{1}{\ell!}\Biggl(\sum_{m=3}^{2N+2}n^{-(m/2)+1}\frac{1}{m!}E^{\nu_0}\Biggl[\Biggl(\sum_{k=1}^{\infty}\sqrt{a_k}e_k(X_1)e_k(Y)\Biggr)^m\Biggr]\\
+n\sum_{j=2}^{N+1}\frac{(-1)^{j-1}}{j}\\
\times\Biggl(\sum_{m=2}^{N+1}n^{-m/2}\frac{1}{m!}E^{\nu_0}\\
\times\Biggl[\Biggl(\sum_{k=1}^{\infty}\sqrt{a_k}e_k(X_1)e_k(Y)\Biggr)^m\Biggr]\Biggr)^j\Biggr)^\ell\Biggr\},A\Biggr]\\
\to 0\qquad as\ n\to\infty,
\end{aligned}
$$

*where the set $A$ on the left-hand side above is either the whole set $H_1$ or the set $\{\|\hat{Y^M}\|_{H_1}\le\sqrt{n}\varepsilon\delta\}\cap B_{n,1/2}$.*

Proof. All has already been proven except for the final assertion with $A=H_1$, which is trivial since the expectation on the set $(\{\|\hat{Y^M}\|_{H_1}\le\sqrt{n}\varepsilon\delta\}\cap B_{n,1/2})^c$ converges to 0 faster than any polynomial order as $n\to\infty$ for any fixed $M\in\mathbf{N}$, which comes from the fact that the integrand is in $L^p$ for some $p>1$ for any fixed $N\in\mathbf{N}$, and $P(\{\|\hat{Y^M}\|_{H_1}\ge\sqrt{n}\varepsilon\delta\}\cup B_{n,1/2}^c)\to 0$ faster than any polynomial order as $n\to\infty$, for any $\delta>0$ and any $\varepsilon>0$. □

Remark 3.1. The assertion of Theorem 3.14 with $A=H_1$ gives us a pure polynomial expansion of $E^{\nu_0^{\otimes\infty}}[\exp(\frac{n}{2}\Psi_2(\frac{S_n}{n},\frac{S_n}{n})),\|\frac{S_n}{n}\|_B<\varepsilon]$, but the summation does not converge as $N\to\infty$. On the other hand, when $A=\{\|Y\|_{H_1}\le\sqrt{n}\varepsilon\delta\}\cap B_{n,1/2}$, the term of the expectation on the left-hand side (if one writes it as a summation with respect to $k_1$ and $k_2$) converges as $N\to\infty$, but instead of giving us a pure polynomial expansion, it only gives us an approximation, for any fixed $N\in\mathbf{N}$, which can be written as a polynomial plus a term which converges to 0 faster than any polynomial order as $n\to\infty$, and a remainder term.

EXAMPLE 3.15. For example, when $N=1$, since

$$E^Y\left[\exp(\tfrac{1}{2}\Psi_2(Y,Y))E^{\nu_0}\left[\left(\sum_{k=1}^{\infty}\sqrt{a_k}e_k(X_1)e_k(Y)\right)^3\right]\right]=0$$

and

$$\det(I_H-\Psi_2)^{-1/2}=E^Y\left[\exp\left(\tfrac{1}{2}\sum_{k=1}^{\infty}a_ke_k(Y)^2\right)\right]$$

as we claimed before, Theorem 3.14 gives us that there exists an $\varepsilon_0>0$ such that for any $\varepsilon\in(0,\varepsilon_0]$,

$$\begin{aligned}
&n\left(E^{\nu_0^{\otimes\infty}}\left[\exp\left(\frac{n}{2}\Psi_2\left(\frac{S_n}{n},\frac{S_n}{n}\right)\right),\left\|\frac{S_n}{n}\right\|_B<\varepsilon\right]-\det(I_H-\Psi_2)^{-1/2}\right)\\
&\quad\to E^Y\left[\exp\left(\frac{1}{2}\sum_{k=1}^{\infty}a_ke_k(Y)^2\right)\right.\\
&\qquad\qquad\times\left(-\frac{1}{8}\left(\sum_{k=1}^{\infty}a_ke_k(Y)^2\right)^2\right.\\
&\qquad\qquad\qquad+\frac{1}{2}E^{\nu_0}\left[\frac{1}{3!}\left(\sum_{k=1}^{\infty}\sqrt{a_k}e_k(X_1)e_k(Y)\right)^3\right]^2\\
&\qquad\qquad\qquad\left.\left.+\frac{1}{4!}E^{\nu_0}\left[\left(\sum_{k=1}^{\infty}\sqrt{a_k}e_k(X_1)e_k(Y)\right)^4\right]\right)\right],
\end{aligned}$$

as $n\to\infty$.

**4. Higher orders.** The study of (2.3) is easy. First, as we observed in Section 1, $\frac{S_n}{\sqrt{n}}\to Y$ in law, where $Y$ is the random variable defined there (or in Section 3). Next, since $\Phi$ is four times continuously Fréchet differentiable by our assumption, we get that $\Psi_3$ and $\Psi_4$ are bounded, so there exist functions $K_3,K_4:B\times B\to\mathbf{R}$ which are bilinear, symmetric and bounded such that $|\Psi_3(y,y,y)|\le\|y\|_BK_3(y,y)$ and $|\Psi_4(y,y,y,y)|\le\|y\|_B^2K_4(y,y)$ for any $y\in B$. Also, by our Assumption A5, there exist constants $\delta_0>0$ and $C_5>0$ such that $|R_5(x^*,y)|\le\|y\|_B^3K_5(y,y)$ for any $y\in B$ with $\|y\|_B<\delta_0$. Now, by using the fact that $|e^x-1|\le|x|e^{|x|}$ and $|e^x-1-x|\le|x|^2e^{|x|}$ for any $x\in\mathbf{R}$, we get that on the set $\{\|\frac{S_n}{n}\|_B<\varepsilon\}$,

$$\left|n\left(\exp\left(\frac{n}{4!}\Psi_4\left(\left(\frac{S_n}{n}\right)^4+nR_5\left(x^*,\frac{S_n}{n}\right)\right)\right)-1\right)-\frac{1}{4!}\Psi_4\left(\left(\frac{S_n}{\sqrt{n}}\right)^4\right)\right|$$

$$
\begin{aligned}
&= \bigg| n\bigg(\exp\bigg(\frac{1}{4!}\frac{1}{n}\Psi_4\bigg(\bigg(\frac{S_n}{\sqrt{n}}\bigg)^4\bigg)\bigg)-1\bigg) - \frac{1}{4!}\Psi_4\bigg(\bigg(\frac{S_n}{\sqrt{n}}\bigg)^4\bigg)\\
&\qquad + \exp\bigg(\frac{1}{4!}\frac{1}{n}\Psi_4\bigg(\bigg(\frac{S_n}{\sqrt{n}}\bigg)^4\bigg)\bigg) n\bigg(\exp\bigg(nR_5\bigg(x^*,\frac{S_n}{n}\bigg)\bigg)-1\bigg)\bigg|\\
&\le \frac{1}{(4!)^2}\frac{1}{n}\bigg\|\frac{S_n}{\sqrt{n}}\bigg\|_B^4 K_4\bigg(\frac{S_n}{\sqrt{n}},\frac{S_n}{\sqrt{n}}\bigg)^2 \exp\bigg(\frac{1}{4!}\varepsilon^2 K_4\bigg(\frac{S_n}{\sqrt{n}},\frac{S_n}{\sqrt{n}}\bigg)\bigg)\\
&\qquad + \exp\bigg(\frac{1}{4!}\varepsilon^2 K_4\bigg(\frac{S_n}{\sqrt{n}},\frac{S_n}{\sqrt{n}}\bigg)\bigg)\exp\bigg(\varepsilon^3 K_5\bigg(\frac{S_n}{\sqrt{n}},\frac{S_n}{\sqrt{n}}\bigg)\bigg)\\
&\qquad\quad \times n^{-1/2}\bigg\|\frac{S_n}{\sqrt{n}}\bigg\|_B^3 K_5\bigg(\frac{S_n}{\sqrt{n}},\frac{S_n}{\sqrt{n}}\bigg),
\end{aligned}
$$

which converges to 0 as $n\to\infty$, $\nu_0^{\otimes\infty}$-a.s. Moreover, by Proposition 3.2, for any $p>1$, there exists an $\varepsilon_0>0$ such that for any $\varepsilon\in(0,\varepsilon_0]$, the right-hand side above is $L^p(\nu_0^{\otimes\infty})$-bounded for $n\in\mathbf{N}$.

Also,

$$
\exp\bigg(\frac{n}{3!}\Psi_3\bigg(\frac{S_n}{n},\frac{S_n}{n},\frac{S_n}{n}\bigg)\bigg)\to 1 \qquad \text{as } n\to\infty,\ \nu_0^{\otimes\infty}\text{-a.s.},
$$

and for any $p>1$, there exists an $\varepsilon_0>0$ such that for any $\varepsilon\in(0,\varepsilon_0]$, it is $L^p$-bounded on the set $\{\|\frac{S_n}{n}\|_B\le\varepsilon\}$ for $n\in\mathbf{N}$.

Therefore, we have the following fact concerning (2.3): there exists a constant $\varepsilon_0>0$ such that for any $\varepsilon\in(0,\varepsilon_0]$,

$$
\begin{aligned}
&E^{\nu_0^{\otimes\infty}}\bigg[\exp\bigg(\frac{n}{2}\Psi_2\bigg(\frac{S_n}{n},\frac{S_n}{n}\bigg)+\frac{n}{3!}\Psi_3\bigg(\frac{S_n}{n},\frac{S_n}{n},\frac{S_n}{n}\bigg)\bigg)\\
&\qquad \times n\bigg(\exp\bigg(\frac{n}{4!}\Psi_4\bigg(\frac{S_n}{n},\frac{S_n}{n},\frac{S_n}{n},\frac{S_n}{n}\bigg)+nR_5\bigg(x^*,\frac{S_n}{n}\bigg)\bigg)-1\bigg),\bigg\|\frac{S_n}{n}\bigg\|_B<\varepsilon\bigg]\\
&\quad\to \frac{1}{4!}E[e^{1/2\Psi_2(Y,Y)}\Psi_4(Y,Y,Y,Y)] \qquad \text{as } n\to\infty.
\end{aligned}
$$

## 5. Third order.

We deal with the term (2.2) in this section.

First, as claimed in Section 4, since $\Psi_3$ is a bounded operator, there exists a function $K_3:B\times B\to\mathbf{R}$ which is bilinear, symmetric and bounded such that $|\Psi_3(y,y,y)|\le\|y\|_B K_3(y,y)$ for any $y\in B$. So by using the fact that $|e^x-1-x-\frac{x^2}{2}|\le e^{|x|}\frac{|x|^3}{3!}$, it is easy to see that

$$
\begin{aligned}
&nE^{\nu_0^{\otimes\infty}}\bigg[\exp\bigg(\frac{1}{2}\Psi_2\bigg(\frac{S_n}{\sqrt{n}},\frac{S_n}{\sqrt{n}}\bigg)\bigg)\\
&\qquad \times\bigg(\exp\bigg(\frac{1}{3!}\frac{1}{\sqrt{n}}\Psi_3\bigg(\bigg(\frac{S_n}{\sqrt{n}}\bigg)^3\bigg)\bigg)-1
\end{aligned}
$$

$$
\begin{aligned}
&- \frac{1}{3!}\frac{1}{\sqrt{n}}\Psi_3\biggl(\biggl(\frac{S_n}{\sqrt{n}}\biggr)^3\biggr) - \frac{1}{2}\biggl(\frac{1}{3!}\frac{1}{\sqrt{n}}\Psi_3\biggl(\biggl(\frac{S_n}{\sqrt{n}}\biggr)^3\biggr)\biggr)^2\biggr), \biggl\|\frac{S_n}{n}\biggr\|_B < \varepsilon\biggr] \\
&\le n^{-1/2} E^{\nu_0^{\otimes\infty}}\biggl[\exp\biggl(\frac{1}{2}\Psi_2\biggl(\frac{S_n}{\sqrt{n}}, \frac{S_n}{\sqrt{n}}\biggr)\biggr)\frac{1}{3!}\exp\biggl(\frac{1}{3!}\varepsilon K_3\biggl(\frac{S_n}{\sqrt{n}}, \frac{S_n}{\sqrt{n}}\biggr)\biggr) \\
&\qquad \times \frac{1}{(3!)^3}\biggl\|\frac{S_n}{\sqrt{n}}\biggr\|_B^3 K_3\biggl(\frac{S_n}{\sqrt{n}}, \frac{S_n}{\sqrt{n}}\biggr)^3, \biggl\|\frac{S_n}{n}\biggr\|_B < \varepsilon\biggr],
\end{aligned}
$$

and the expectation on the right-hand side above is bounded for $n \in \mathbf{N}$ if $\varepsilon > 0$ is small enough, by Proposition 3.2. Therefore, we get that there exists an $\varepsilon_0 > 0$ such that for any $\varepsilon \in (0, \varepsilon_0]$, the left-hand side above converges to 0 as $n \to \infty$.

Also, as observed in Section 1, $\frac{S_n}{\sqrt{n}} \to Y$ in law as $n \to \infty$, where $Y$ is the random variable defined there (or in Section 3). Therefore, there exists an $\varepsilon_0 > 0$ such that for any $\varepsilon \in (0, \varepsilon_0]$, as $n \to \infty$,

$$
\begin{aligned}
&n E^{\nu_0^{\otimes\infty}}\biggl[\exp\biggl(\frac{1}{2}\Psi_2\biggl(\frac{S_n}{\sqrt{n}}, \frac{S_n}{\sqrt{n}}\biggr)\biggr)\frac{1}{2}\biggl(\frac{1}{3!}\frac{1}{\sqrt{n}}\Psi_3\biggl(\biggl(\frac{S_n}{\sqrt{n}}\biggr)^3\biggr)\biggr)^2, \biggl\|\frac{S_n}{n}\biggr\|_B < \varepsilon\biggr] \\
&\quad = E^{\nu_0^{\otimes\infty}}\biggl[\exp\biggl(\frac{1}{2}\Psi_2\biggl(\frac{S_n}{\sqrt{n}}, \frac{S_n}{\sqrt{n}}\biggr)\biggr)\frac{1}{2(3!)^2}\Psi_3\biggl(\biggl(\frac{S_n}{\sqrt{n}}\biggr)^3\biggr)^2, \biggl\|\frac{S_n}{n}\biggr\|_B < \varepsilon\biggr] \\
&\quad \to \frac{1}{2(3!)^2} E[e^{\Psi_2(Y,Y)/2}\Psi_3(Y,Y,Y)^2],
\end{aligned}
$$

where in the last line, we used Proposition 3.1 and the general fact that a.s.-convergence and $L^p$-boundedness for some $p > 1$ imply $L^1$-convergence.

Therefore, we only need to study the term

$$
\begin{aligned}
&\frac{1}{3!}\sqrt{n} E^{\nu_0^{\otimes\infty}}\biggl[\exp\biggl(\frac{1}{2}\Psi_2\biggl(\frac{S_n}{\sqrt{n}}, \frac{S_n}{\sqrt{n}}\biggr)\biggr)\Psi_3\biggl(\biggl(\frac{S_n}{\sqrt{n}}\biggr)^3\biggr), \biggl\|\frac{S_n}{n}\biggr\|_B < \varepsilon\biggr] \\
&\quad = \frac{1}{3!}\sqrt{n} E^{\nu_0^{\otimes\infty}}\biggl[E^Y\biggl[\exp\biggl(\sum_{k=1}^{\infty}\sqrt{a_k} e_k\biggl(\frac{S_n}{\sqrt{n}}\biggr) e_k(Y)\biggr)\biggr]\Psi_3\biggl(\biggl(\frac{S_n}{\sqrt{n}}\biggr)^3\biggr), \\
&\qquad\qquad \biggl\|\frac{S_n}{n}\biggr\|_B < \varepsilon\biggr],
\end{aligned}
$$

for $\varepsilon > 0$ small enough, where $a_k$, $e_k$, $k \in \mathbf{N}$, and $Y$ are as defined in Section 1.

First, and similarly as in Section 3, we have the following three propositions (with the same notation as there).

PROPOSITION 5.1. *There exists a $\varepsilon_0 > 0$ such that for any $\varepsilon \in (0, \varepsilon_0]$ and any $\varepsilon_1 > 0$,*

$$
\sqrt{n} E^{\nu_0^{\otimes\infty}}\biggl[\exp\biggl(\frac{1}{2}\Psi_2\biggl(\frac{S_n}{\sqrt{n}}, \frac{S_n}{\sqrt{n}}\biggr)\biggr)\Psi_3\biggl(\frac{S_n}{\sqrt{n}}, \frac{S_n}{\sqrt{n}}, \frac{S_n}{\sqrt{n}}\biggr),
$$

$$\left\{\left\|\frac{S_n}{n}\right\|_B > \varepsilon_1\right\} \cap \left\{\left\|\frac{S_n}{n}\right\|_{H_M} < \varepsilon\right\}\Bigg]$$

*converges to* 0 *exponentially as* $n \to \infty$.

PROOF. The proof goes in the same way as the one of Proposition 3.3. (We remark again that $\Psi_3$ is bounded by our assumption.) □

PROPOSITION 5.2. *There exists an* $\varepsilon_0 \geq 0$ *such that for any* $\varepsilon \in (0, \varepsilon_0]$ *and any* $\delta > 0$,

$$\sqrt{n} E^Y\left[E^{\nu_0^{\otimes\infty}}\left[\exp\left(\sum_{k=1}^{\infty}\sqrt{a_k} e_k\left(\frac{S_n}{\sqrt{n}}\right) e_k(Y)\right)\Psi_3\left(\left(\frac{S_n}{\sqrt{n}}\right)^3\right), \left\|\frac{S_n}{n}\right\|_{H_M} < \varepsilon\right],\right.$$
$$\left.\|Y^{\hat{M}}\|_{H_1} \geq \sqrt{n}\varepsilon\delta\right]$$

*converges to* 0 *exponentially as* $n \to \infty$.

PROOF. It suffices to remark that $n \leq \frac{\|Y^{\hat{M}}\|_{H_1}^2}{\varepsilon^2\delta^2}$ on the set $\{\|Y^{\hat{M}}\|_{H_1} \geq \sqrt{n}\varepsilon\delta\}$. The proof goes then in the same way as the one of Proposition 3.3. □

PROPOSITION 5.3. *There exists a* $\delta_0 \geq 0$ *such that for any* $\delta \in (0, \delta_0]$ *and any* $\varepsilon > 0$,

$$\sqrt{n} E^Y\left[E^{\nu_0^{\otimes\infty}}\left[\exp\left(\sum_{k=1}^{\infty}\sqrt{a_k} e_k\left(\frac{S_n}{\sqrt{n}}\right) e_k(Y)\right)\Psi_3\left(\left(\frac{S_n}{\sqrt{n}}\right)^3\right), \left\|\frac{S_n}{n}\right\|_{H_M} \geq \varepsilon\right],\right.$$
$$\left.\|Y^{\hat{M}}\|_{H_1} < \sqrt{n}\varepsilon\delta\right]$$

*converges to* 0 *exponentially as* $n \to \infty$.

PROOF. The idea is the same as before, except that this time, we have $n \leq \frac{1}{\varepsilon^2}\|\frac{S_n}{\sqrt{n}}\|_{H_M}^2$ on the set $\{\|\frac{S_n}{n}\|_{H_M} \geq \varepsilon\}$. □

By the above results, we only need to discuss the asymptotic expansion in powers of $\frac{1}{n^{1/2}}$ as $n \to \infty$ of

$$\frac{1}{3!}\sqrt{n} E^Y\left[E^{\nu_0^{\otimes\infty}}\left[\exp\left(\sum_{k=1}^{\infty}\sqrt{a_k} e_k\left(\frac{S_n}{\sqrt{n}}\right) e_k(Y)\right)\Psi_3\left(\left(\frac{S_n}{\sqrt{n}}\right)^3\right)\right],\right.$$
$$\left.\|Y^{\hat{M}}\|_{H_1} < \sqrt{n}\varepsilon\delta\right]$$

for $\varepsilon > 0$ and $\delta > 0$ small enough.

In the same way as for Proposition 3.8, we have the following.

PROPOSITION 5.4. *For any* $\xi > 0$,

$$E^Y\left[E^{\nu_0}\left[\exp\left(\frac{1}{\sqrt{n}}\sum_{k=1}^{\infty}\sqrt{a_k}e_k(X_1)e_k(Y)\right)\right]^n,\right.$$

$$\left.\{\|Y^{\hat{M}}\|_{H_1} < \sqrt{n}\varepsilon\delta\} \cap B_{n,\xi}^c\right] \to 0$$

*as* $n \to \infty$, *faster than any polynomial order.* □

Therefore, we only need to discuss the asymptotic expansion of

$$\frac{1}{3!}\sqrt{n}E^Y\left[E^{\nu_0^{\otimes\infty}}\left[\exp\left(\sum_{k=1}^{\infty}\sqrt{a_k}e_k\left(\frac{S_n}{\sqrt{n}}\right)e_k(Y)\right)\Psi_3\left(\left(\frac{S_n}{\sqrt{n}}\right)^3\right)\right],\right.$$

$$\left.\{\|Y^{\hat{M}}\|_{H_1} < \sqrt{n}\varepsilon\delta\} \cap B_{n,1/2}\right]$$

for $\varepsilon > 0$ and $\delta > 0$ small enough. We notice that this expectation can be decomposed as follows:

$$\frac{1}{3!}\frac{1}{n}E^Y\left[E^{\nu_0^{\otimes\infty}}\left[\exp\left(\sum_{k=1}^{\infty}\sqrt{a_k}e_k\left(\frac{S_n}{\sqrt{n}}\right)e_k(Y)\right)\Psi_3(S_n,S_n,S_n)\right],\right.$$

$$\left.\{\|Y^{\hat{M}}\|_{H_1} < \sqrt{n}\varepsilon\delta\} \cap B_{n,1/2}\right]$$

$$= \frac{1}{3!}E^Y\left[E^{\nu_0^{\otimes\infty}}\left[\exp\left(\sum_{k=1}^{\infty}\sqrt{a_k}e_k\left(\frac{S_{n-1}}{\sqrt{n}}\right)e_k(Y)\right)\right]\right.$$

$$\times E^{\nu_0}\left[\exp\left(\sum_{k=1}^{\infty}\sqrt{a_k}e_k\left(\frac{X_1}{\sqrt{n}}\right)e_k(Y)\right)\Psi_3(X_1,X_1,X_1)\right],$$

$$\left.\{\|Y^{\hat{M}}\|_{H_1} < \sqrt{n}\varepsilon\delta\} \cap B_{n,1/2}\right]$$

$$+ \frac{1}{2}(n-1)E^Y\left[E^{\nu_0^{\otimes\infty}}\left[\exp\left(\sum_{k=1}^{\infty}\sqrt{a_k}e_k\left(\frac{S_{n-2}}{\sqrt{n}}\right)e_k(Y)\right)\right]\right.$$

$$\times E^{\nu_0^{\otimes 2}}\left[\exp\left(\sum_{i=1}^{2}\sum_{k=1}^{\infty}\sqrt{a_k}e_k\left(\frac{X_i}{\sqrt{n}}\right)e_k(Y)\right)\right. \tag{5.1}$$

$$\left.\left.\times \Psi_3(X_1, X_2, X_2)\right]\right],$$

$$\left.\{\|Y^{\hat{M}}\|_{H_1} < \sqrt{n}\varepsilon\delta\} \cap B_{n,1/2}\right]$$

$$+ \frac{1}{3!}(n-1)(n-2)$$

$$\times E^Y\left[E^{\nu_0^{\otimes\infty}}\left[\exp\left(\sum_{k=1}^{\infty}\sqrt{a_k}e_k\left(\frac{S_{n-3}}{\sqrt{n}}\right)e_k(Y)\right)\right]\right.$$

$$\left.\times E^{\nu_0^{\otimes 3}}\left[\exp\left(\sum_{i=1}^{3}\sum_{k=1}^{\infty}\sqrt{a_k}e_k\left(\frac{X_i}{\sqrt{n}}\right)e_k(Y)\right)\Psi_3(X_1, X_2, X_3)\right],\right.$$

$$\left.\{\|Y^{\hat{M}}\|_{H_1} < \sqrt{n}\varepsilon\delta\} \cap B_{n,1/2}\right].$$

From now on, we study each of the three terms above. As before, we show the convergence for a.e.-$Y$; then with the fact that there exists a $p > 1$ such that the three integrands in (5.1) are $L^p$ (with respect to the distribution of $Y$)-bounded on the set $\{\|Y^{\hat{M}}\|_{H_1} < \sqrt{n}\varepsilon\delta\} \cap B_{n,1/2}$ with respect to $n \in \mathbf{N}$, and the dominated convergence theorem, we get our assertion. The proof of the $L^p$-boundedness (for some $p > 1$) is similar to the one given before, and therefore we will omit it.

It is easy to see by Theorem 3.14 that for a.e.-$Y$,

$$\sqrt{n}\left(E^{\nu_0}\left[\exp\left(\sum_{k=1}^{\infty}\sqrt{a_k}e_k\left(\frac{X_1}{\sqrt{n}}\right)e_k(Y)\right)\right]^n - \exp\left(\frac{1}{2}\sum_{k=1}^{\infty}a_k e_k(Y)^2\right)\right)$$

$$\to \exp\left(\frac{1}{2}\sum_{k=1}^{\infty}a_k e_k(Y)^2\right)\frac{1}{3!}E^{\nu_0}\left[\left(\sum_{k=1}^{\infty}\sqrt{a_k}e_k(X_1)e_k(Y)\right)^3\right]$$

as $n \to \infty$. Hence, for $i = 1, 2, 3$,

$$E^{\nu_0}\left[\exp\left(\sum_{k=1}^{\infty}\sqrt{a_k}e_k\left(\frac{X_1}{\sqrt{n}}\right)e_k(Y)\right)\right]^{n-i} \to \exp\left(\frac{1}{2}\sum_{k=1}^{\infty}a_k e_k(Y)^2\right) \tag{5.2}$$

and

$$\sqrt{n}\left(1 - E^{\nu_0}\left[\exp\left(\sum_{k=1}^{\infty}\sqrt{a_k}e_k\left(\frac{X_1}{\sqrt{n}}\right)e_k(Y)\right)\right]^i\right) \to 0.$$

Therefore,

$$\sqrt{n}\left(E^{\nu_0^{\otimes\infty}}\left[\exp\left(\sum_{k=1}^{\infty}\sqrt{a_k}e_k\left(\frac{S_{n-i}}{\sqrt{n}}\right)e_k(Y)\right)\right] - \exp\left(\frac{1}{2}\sum_{k=1}^{\infty}a_k e_k(Y)^2\right)\right)$$

$$= \sqrt{n}\left(E^{\nu_0}\left[\exp\left(\sum_{k=1}^{\infty}\sqrt{a_k}e_k\left(\frac{X_1}{\sqrt{n}}\right)e_k(Y)\right)\right]^n - \exp\left(\frac{1}{2}\sum_{k=1}^{\infty}a_k e_k(Y)^2\right)\right)$$

$$+ E^{\nu_0}\left[\exp\left(\sum_{k=1}^{\infty}\sqrt{a_k}e_k\left(\frac{X_1}{\sqrt{n}}\right)e_k(Y)\right)\right]^{n-i}$$

$$\times \sqrt{n}\left(1 - E^{\nu_0}\left[\exp\left(\sum_{k=1}^{\infty}\sqrt{a_k}e_k\left(\frac{X_1}{\sqrt{n}}\right)e_k(Y)\right)\right]^i\right)$$

$$\to \frac{1}{3!}\exp\left(\frac{1}{2}\sum_{k=1}^{\infty}a_k e_k(Y)^2\right)E^{\nu_0}\left[\left(\sum_{k=1}^{\infty}\sqrt{a_k}e_k(X_1)e_k(Y)\right)^3\right],$$

as $n \to \infty$.

The first term in the decomposition (5.1) converges to

$$\frac{1}{3!}E^Y\left[\exp\left(\frac{1}{2}\sum_{k=1}^{\infty}a_k e_k(Y)^2\right)\right]E^{\nu_0}[\Psi_3(X_1,X_1,X_1)],$$

by (5.2) and the easy fact that

$$E^{\nu_0}\left[\exp\left(\sum_{k=1}^{\infty}\sqrt{a_k}e_k\left(\frac{X_1}{\sqrt{n}}\right)e_k(Y)\right)\Psi_3(X_1,X_1,X_1)\right]$$
$$\to E^{\nu_0}[\Psi_3(X_1,X_1,X_1)] \qquad \text{as } n \to \infty.$$

For the second term in the decomposition (5.1), we notice that since $\Psi_3$ is trilinear and $\nu_0$ has mean 0, we have that $E^{\nu_0}[\Psi_3(X_1,\cdot,\cdot)] = 0$, so

$$\sqrt{n}E^{\nu_0^{\otimes 2}}\left[\exp\left(\frac{1}{\sqrt{n}}\sum_{i=1}^{2}\sum_{k=1}^{\infty}\sqrt{a_k}e_k(X_i)e_k(Y)\right)\Psi_3(X_1,X_2,X_2)\right]$$
$$\to E^{\nu_0^{\otimes 2}}\left[\left(\sum_{k=1}^{\infty}\sqrt{a_k}e_k(X_1)e_k(Y)\right)\Psi_3(X_1,X_2,X_2)\right] \qquad \text{as } n \to \infty.$$

Moreover, since

$$E^Y\left[\exp\left(\tfrac{1}{2}\sum_{k=1}^{\infty}a_k e_k(Y)^2\right)e_j(Y)\right] = 0 \qquad \text{for } j \in \mathbf{N},$$

we have also that

$$nE^Y\left[\exp\left(\frac{1}{2}\sum_{k=1}^{\infty}a_k e_k(Y)^2\right)\right.$$
$$\left.\times E^{\nu_0^{\otimes 2}}\left[\exp\left(\frac{1}{\sqrt{n}}\sum_{i=1}^{2}\sum_{k=1}^{\infty}\sqrt{a_k}e_k(X_i)e_k(Y)\right)\Psi_3(X_1,X_2,X_2)\right]\right.,$$

$$\left.\{\|Y^{\hat{M}}\|_{H_1} < \sqrt{n}\varepsilon\delta\} \cap B_{n,1/2}\right]$$

$$\begin{aligned}
&\to E^Y\Bigg[\exp\bigg(\frac{1}{2}\sum_{k=1}^{\infty} a_k e_k(Y)^2\bigg) \\
&\quad\times E^{\nu_0^{\otimes 2}}\Bigg[\bigg(\sum_{k=1}^{\infty}\sqrt{a_k}e_k(X_1)e_k(Y)\bigg) \\
&\qquad\times \bigg(\sum_{k=1}^{\infty}\sqrt{a_k}e_k(X_2)e_k(Y)\bigg)\Psi_3(X_1,X_2,X_2)\Bigg]\Bigg] \\
&\quad+\frac{1}{2}E^Y\Bigg[\exp\bigg(\frac{1}{2}E^{\nu_0}\bigg[\bigg(\sum_{k=1}^{\infty}\sqrt{a_k}e_k(X_1)e_k(Y)\bigg)^2\bigg]\bigg) \\
&\qquad\times E^{\nu_0^{\otimes 2}}\Bigg[\bigg(\sum_{k=1}^{\infty}\sqrt{a_k}e_k(X_1)e_k(Y)\bigg)^2\Psi_3(X_1,X_2,X_2)\Bigg]\Bigg].
\end{aligned}$$

Therefore, we have that the second term of (5.1) is equal

$$\begin{aligned}
&\frac{1}{2}\frac{n-1}{n}E^Y\Bigg[\sqrt{n}\Bigg(E^{\nu_0^{\otimes\infty}}\bigg[\exp\bigg(\frac{1}{\sqrt{n}}\sum_{k=1}^{\infty}\sqrt{a_k}e_k(X_1)e_k(Y)\bigg)\bigg]^{n-2} \\
&\qquad\qquad -\exp\bigg(\frac{1}{2}\sum_{k=1}^{\infty}a_k e_k(Y)^2\bigg)\Bigg) \\
&\qquad\times\sqrt{n}E^{\nu_0^{\otimes 2}}\bigg[\exp\bigg(\frac{1}{\sqrt{n}}\sum_{i=1}^{2}\sum_{k=1}^{\infty}\sqrt{a_k}e_k(X_i)e_k(Y)\bigg)\Psi_3(X_1,X_2,X_2)\bigg], \\
&\qquad\qquad \left.\{\|Y^{\hat{M}}\|_{H_1} < \sqrt{n}\varepsilon\delta\}\cap B_{n,1/2}\right] \\
&\quad+\frac{1}{2}\frac{n-1}{n}nE^Y\Bigg[\exp\bigg(\frac{1}{2}\sum_{k=1}^{\infty}a_k e_k(Y)^2\bigg) \\
&\qquad\times E^{\nu_0^{\otimes 2}}\bigg[\exp\bigg(\frac{1}{\sqrt{n}}\sum_{i=1}^{2}\sum_{k=1}^{\infty}\sqrt{a_k}e_k(X_i)e_k(Y)\bigg)\Psi_3(X_1,X_2,X_2)\bigg], \\
&\qquad\qquad \left.\{\|Y^{\hat{M}}\|_{H_1} < \sqrt{n}\varepsilon\delta\}\cap B_{n,1/2}\right] \\
&\to\frac{1}{2\cdot 3!}E^Y\Bigg[\exp\bigg(\frac{1}{2}\sum_{k=1}^{\infty}a_k e_k(Y)^2\bigg)E^{\nu_0}\bigg[\bigg(\sum_{k=1}^{\infty}\sqrt{a_k}e_k(X_1)e_k(Y)\bigg)^3\bigg]
\end{aligned}$$

$$\times E^{\nu_0^{\otimes 2}}\left[\left(\sum_{k=1}^{\infty}\sqrt{a_k}e_k(X_1)e_k(Y)\right)\Psi_3(X_1,X_2,X_2)\right]\Bigg]$$

$$+\frac{1}{2}E^Y\Bigg[\exp\left(\frac{1}{2}\sum_{k=1}^{\infty}a_ke_k(Y)^2\right)$$

$$\times E^{\nu_0^{\otimes 2}}\Bigg[\left(\sum_{k=1}^{\infty}\sqrt{a_k}e_k(X_1)e_k(Y)\right)$$

$$\times\left(\sum_{k=1}^{\infty}\sqrt{a_k}e_k(X_2)e_k(Y)\right)\Psi_3(X_1,X_2,X_2)\Bigg]\Bigg]$$

$$+\frac{1}{4}E^Y\Bigg[\exp\left(\frac{1}{2}\sum_{k=1}^{\infty}a_ke_k(Y)^2\right)$$

$$\times E^{\nu_0^{\otimes 2}}\left[\left(\sum_{k=1}^{\infty}\sqrt{a_k}e_k(X_1)e_k(Y)\right)^2\Psi_3(X_1,X_2,X_2)\right]\Bigg].$$

The third term of (5.1) can be dealt with in the same way. We have that

$$n^{3/2}E^{\nu_0^{\otimes 3}}\left[\exp\left(\frac{1}{\sqrt{n}}\sum_{i=1}^{3}\sum_{k=1}^{\infty}\sqrt{a_k}e_k(X_i)e_k(Y)\right)\Psi_3(X_1,X_2,X_3)\right]$$

$$\to E^{\nu_0^{\otimes 3}}\left[\prod_{i=1}^{3}\left(\sum_{k=1}^{\infty}\sqrt{a_k}e_k(X_i)e_k(Y)\right)\Psi_3(X_1,X_2,X_3)\right]$$

and

$$n^2E^Y\Bigg[\exp\left(\frac{1}{2}\sum_{k=1}^{\infty}a_ke_k(Y)^2\right)$$

$$\times E^{\nu_0^{\otimes 3}}\left[\exp\left(\frac{1}{\sqrt{n}}\sum_{i=1}^{3}\sum_{k=1}^{\infty}\sqrt{a_k}e_k(X_i)e_k(Y)\right)\Psi_3(X_1,X_2,X_3)\right],$$

$$\{\|Y^{\hat{M}}\|_{H_1}<\sqrt{n}\varepsilon\delta\}\cap B_{n,1/2}\Bigg]$$

$$\to 3E^Y\Bigg[\exp\left(\frac{1}{2}\sum_{k=1}^{\infty}a_ke_k(Y)^2\right)$$

$$\times E^{\nu_0^{\otimes 3}}\Bigg[\left(\sum_{k=1}^{\infty}\sqrt{a_k}e_k(X_1)e_k(Y)\right)\left(\sum_{k=1}^{\infty}\sqrt{a_k}e_k(X_2)e_k(Y)\right)$$

$$\times\left(\sum_{k=1}^{\infty}\sqrt{a_k}e_k(X_3)e_k(Y)\right)^2\Psi_3(X_1,X_2,X_3)\Bigg]\Bigg].$$

Therefore, the third term of (5.1) is equal

$$
\begin{aligned}
&\frac{1}{3!}\frac{(n-1)(n-2)}{n^2}\\
&\quad\times E^Y\Bigg[\sqrt{n}\Bigg(E^{\nu_0^{\otimes\infty}}\Bigg[\exp\Bigg(\frac{1}{\sqrt{n}}\sum_{k=1}^{\infty}\sqrt{a_k}e_k(X_1)e_k(Y)\Bigg)\Bigg]^{n-3}\\
&\qquad\qquad\qquad\qquad\qquad -\exp\Bigg(\frac{1}{2}\sum_{k=1}^{\infty}a_k e_k(Y)^2\Bigg)\Bigg)\\
&\qquad\qquad \times n^{3/2}E^{\nu_0^{\otimes 3}}\Bigg[\exp\Bigg(\frac{1}{\sqrt{n}}\sum_{i=1}^{3}\sum_{k=1}^{\infty}\sqrt{a_k}e_k(X_i)e_k(Y)\Bigg)\Psi_3(X_1,X_2,X_3)\Bigg],\\
&\qquad\qquad\qquad\qquad\qquad\qquad \{\|Y^{\hat{M}}\|_{H_1}<\sqrt{n}\varepsilon\delta\}\cap B_{n,1/2}\Bigg]\\
&\quad+\frac{1}{3!}\frac{(n-1)(n-2)}{n^2}\\
&\qquad\times n^2E^Y\Bigg[\exp\Bigg(\frac{1}{2}\sum_{k=1}^{\infty}a_k e_k(Y)^2\Bigg)\\
&\qquad\qquad\times E^{\nu_0^{\otimes 3}}\Bigg[\exp\Bigg(\frac{1}{\sqrt{n}}\sum_{i=1}^{3}\sum_{k=1}^{\infty}\sqrt{a_k}e_k(X_i)e_k(Y)\Bigg)\Psi_3(X_1,X_2,X_3)\Bigg],\\
&\qquad\qquad\qquad\qquad\qquad\qquad \{\|Y^{\hat{M}}\|_{H_1}<\sqrt{n}\varepsilon\delta\}\cap B_{n,1/2}\Bigg]\\
&\quad\to\frac{1}{(3!)^2}E^Y\Bigg[\exp\Bigg(\frac{1}{2}\sum_{k=1}^{\infty}a_k e_k(Y)^2\Bigg)E^{\nu_0}\Bigg[\Bigg(\sum_{k=1}^{\infty}\sqrt{a_k}e_k(X_i)e_k(Y)\Bigg)^3\Bigg]\\
&\qquad\qquad\times E^{\nu_0^{\otimes 3}}\Bigg[\prod_{i=1}^{3}\Bigg(\sum_{k=1}^{\infty}\sqrt{a_k}e_k(X_i)e_k(Y)\Bigg)\Psi_3(X_1,X_2,X_3)\Bigg]\Bigg]\\
&\qquad+\frac{1}{2}E^Y\Bigg[\exp\Bigg(\frac{1}{2}\sum_{k=1}^{\infty}a_k e_k(Y)^2\Bigg)\\
&\qquad\qquad\times E^{\nu_0^{\otimes 3}}\Bigg[\Bigg(\sum_{k=1}^{\infty}\sqrt{a_k}e_k(X_1)e_k(Y)\Bigg)\Bigg(\sum_{k=1}^{\infty}\sqrt{a_k}e_k(X_2)e_k(Y)\Bigg)\\
&\qquad\qquad\qquad\times\Bigg(\sum_{k=1}^{\infty}\sqrt{a_k}e_k(X_3)e_k(Y)\Bigg)^2\Psi_3(X_1,X_2,X_3)\Bigg]\Bigg],
\end{aligned}
$$

as $n\to\infty$.

In conclusion, we have proven the following.

LEMMA 5.5. *There exists an* $\varepsilon_0 > 0$ *such that for any* $\varepsilon \in (0, \varepsilon_0)$,

$$\lim_{n\to\infty} nE^{\nu_0^{\otimes\infty}}\left[\exp\left(\frac{n}{2}\Psi_2\left(\frac{S_n}{n},\frac{S_n}{n}\right)\right)\left(\exp\left(\frac{n}{3!}\Psi_3\left(\frac{S_n}{n},\frac{S_n}{n},\frac{S_n}{n}\right)\right)-1\right), \left\|\frac{S_n}{n}\right\|_B < \varepsilon\right]$$

$$= \frac{1}{2(3!)^2}E[e^{\Psi_2(Y,Y)/2}\Psi_3(Y,Y,Y)^2]$$

$$+ \frac{1}{3!}E^Y\left[\exp\left(\frac{1}{2}\sum_{k=1}^\infty a_k e_k(Y)^2\right)\right]E^{\nu_0}[\Psi_3(X_1,X_1,X_1)]$$

$$+ \frac{1}{2\cdot 3!}E^Y\left[\exp\left(\frac{1}{2}\sum_{k=1}^\infty a_k e_k(Y)^2\right)E^{\nu_0}\left[\left(\sum_{k=1}^\infty \sqrt{a_k}e_k(X_1)e_k(Y)\right)^3\right] \times E^{\nu_0^{\otimes 2}}\left[\left(\sum_{k=1}^\infty \sqrt{a_k}e_k(X_1)e_k(Y)\right)\Psi_3(X_1,X_2,X_2)\right]\right]$$

$$+ \frac{1}{2}E^Y\left[\exp\left(\frac{1}{2}\sum_{k=1}^\infty a_k e_k(Y)^2\right) \times E^{\nu_0^{\otimes 2}}\left[\left(\sum_{k=1}^\infty \sqrt{a_k}e_k(X_1)e_k(Y)\right) \times \left(\sum_{k=1}^\infty \sqrt{a_k}e_k(X_2)e_k(Y)\right)\Psi_3(X_1,X_2,X_2)\right]\right]$$

$$+ \frac{1}{4}E^Y\left[\exp\left(\frac{1}{2}\sum_{k=1}^\infty a_k e_k(Y)^2\right) \times E^{\nu_0^{\otimes 2}}\left[\left(\sum_{k=1}^\infty \sqrt{a_k}e_k(X_1)e_k(Y)\right)^2\Psi_3(X_1,X_2,X_2)\right]\right]$$

$$+ \frac{1}{(3!)^2}E^Y\left[\exp\left(\frac{1}{2}\sum_{k=1}^\infty a_k e_k(Y)^2\right)E^{\nu_0}\left[\left(\sum_{k=1}^\infty \sqrt{a_k}e_k(X_i)e_k(Y)\right)^3\right] \times E^{\nu_0^{\otimes 3}}\left[\prod_{i=1}^3\left(\sum_{k=1}^\infty \sqrt{a_k}e_k(X_i)e_k(Y)\right)\Psi_3(X_1,X_2,X_3)\right]\right]$$

$$+\frac{1}{2}E^Y\left[\exp\left(\frac{1}{2}\sum_{k=1}^{\infty}a_k e_k(Y)^2\right)\right.$$
$$\times E^{\nu_0^{\otimes 3}}\left[\left(\sum_{k=1}^{\infty}\sqrt{a_k}e_k(X_1)e_k(Y)\right)\left(\sum_{k=1}^{\infty}\sqrt{a_k}e_k(X_2)e_k(Y)\right)\right.$$
$$\left.\left.\times\left(\sum_{k=1}^{\infty}\sqrt{a_k}e_k(X_3)e_k(Y)\right)^2\Psi_3(X_1,X_2,X_3)\right]\right].$$

This concludes the proof of Theorem 1.3.

**Acknowledgments.** We are very grateful to Professor D. Sonderman for the financial support given to the second author. The warm hospitality of Professor Luciano Tubaro and the Mathematics Department of the University of Trento, where part of this work was done, is also gratefully acknowledged.

Institute of Applied Mathematics
University of Bonn
Wegelerstrasse 6
D53115 Bonn
Germany
e-mail: albeverio@uni-bonn.de

Graduate School of Information Sciences
Tohoku University
Aramaki aza Aoba
980-8579 Sendai
Japan
e-mail: liang@math.is.tohoku.ac.jp